\documentclass[11pt]{article}
\usepackage[a4paper,margin=27mm]{geometry}
\usepackage{amsmath,amssymb,amsthm,mathtools,bm}
\usepackage{microtype,booktabs,enumitem}
\usepackage{xcolor}
\usepackage[normalem]{ulem}
\usepackage{comment}

\usepackage[numbers,sort&compress]{natbib}
\usepackage[colorlinks=true,linkcolor=blue,citecolor=blue,urlcolor=blue]{hyperref}
\numberwithin{equation}{section}
\newtheorem{theorem}{Theorem}[section]
\newtheorem{proposition}[theorem]{Proposition}
\newtheorem{lemma}[theorem]{Lemma}
\newtheorem{corollary}[theorem]{Corollary}
\theoremstyle{definition}
\newtheorem{definition}[theorem]{Definition}
\newtheorem{assumption}[theorem]{Assumption}
\newtheorem{example}[theorem]{Example}
\theoremstyle{remark}
\newtheorem{remark}[theorem]{Remark}
\newcommand{\E}{\mathbb E}
\newcommand{\R}{\mathbb R}
\newcommand{\F}{\mathcal F}
\newcommand{\Law}{\mathcal L}
\newcommand{\dd}{\mathrm d}
\newcommand{\dt}{\,\mathrm dt}
\newcommand{\ds}{\,\mathrm ds}
\newcommand{\dW}{\,\mathrm dW}
\newcommand{\norm}[1]{\left\|#1\right\|}
\newcommand{\BL}{\mathrm{BL}}
\newcommand{\cP}{\mathcal P}
\newcommand{\eps}{\varepsilon}
\newcommand{\rev}[1]{{\color{red}#1}}

\title{Global Well-Posedness for McKean--Vlasov Control   with Discrete Path Dependence and Convex Law Penalties}
\author{Du OUYANG \footnote{École Polytechnique, CMAP,
\sf du.ouyang at polytechnique.edu.
This author is supported by the Soci\'et\'e G\'en\'erale Chair ``Financial Risks''}  \and Huy\^en PHAM\footnote{École Polytechnique, CMAP,
\sf huyen.pham at polytechnique.edu.
This author is supported by the Soci\'et\'e G\'en\'erale Chair ``Financial Risks'',
by the Chair AIDA-F with Bank of America, by FiME (Laboratory of Finance and Energy Markets) and the EDF–CACIB Chair ``Finance and Sustainable Development''.}}
\date{November 2026}
\begin{document}
\setlength{\emergencystretch}{2em}
\maketitle

\begin{abstract}
Motivated by time-series generation, we study McKean--Vlasov control  with drift, uncontrolled diffusion and running cost depending on current and recorded states at finitely many observation dates, and penalties on their joint laws. Under uniform ellipticity, regularity and Hamiltonian assumptions, convexity of the penalties along mixtures of probability measures yields a unique optimal control in Lipschitz feedback form on every finite horizon. The forward-backward system 
associated with Pontryagin's stochastic maximum principle has adjoint jumps at intermediate observation dates, and admits  a unique square-integrable solution with bounded adjoint $Y$. 
We construct the optimal control via an auxiliary potential
mean field game. Linearizing the penalties at a frozen population
path law gives the agent's problem. We solve its
nonlinear Hamilton-Jacobi-Bellman equations with past observations
as parameters. Uniform Bernstein gradient bounds and mixed derivative
estimates give a continuous best-response map. 
A Schauder-Tychonoff fixed point yields a Nash equilibrium of the game. 
By verification and convexity along mixtures, the equilibrium control minimizes the original McKean-Vlasov objective.    
The framework covers squared maximum mean discrepancy penalties with suitably smooth kernels. With the control as drift, a standard quadratic control cost, current-state diffusion and running marginal penalties, we establish gradient domination and  obtain a posteriori bounds on the cost gap and the squared errors in the control and state in terms of the Pontryagin residual. 
\end{abstract}

\medskip
\noindent\textbf{Key words.} McKean-Vlasov control, discrete path dependence, forward-backward stochastic differential equations with adjoint jumps, potential mean field games, convexity along mixtures, 
Hamilton-Jacobi-Bellman. 

\medskip
\noindent\textbf{MSC codes.} 93E20, 49N80, 60H10, 49L12

\section{Introduction}
Generative models for time series are often trained on trajectories sampled at finitely many dates.
Two stochastic processes can have identical marginal laws at every observation date and still differ in their temporal dependence. This motivates penalizing discrepancies between the joint law of the generated states at those dates and the corresponding target law. 
El Boustany et al.~\cite{elbou2026} incorporate marginal and joint discrepancies in a stochastic control objective using maximum mean discrepancy (MMD) penalties. They use the resulting McKean--Vlasov forward--backward stochastic differential equation (FBSDE) to learn generative dynamics from data. Its adjoint jumps at the intermediate observation dates encode the sensitivities to the recorded states, as characterized by the stochastic maximum principle of~\citet{hafsi2026}.

We study a general McKean–Vlasov control problem of this form.
In the controlled dynamics~\eqref{eq:state}, the drift and the uncontrolled diffusion matrix may depend on the current state and on the observations already recorded. The running cost may depend on the same variables and on the control. The running law penalty acts on the joint law of the current state and the recorded states, and the terminal penalty acts on the joint law of all observations. The penalties are assumed to be convex along mixtures of probability measures, a property satisfied by squared MMD penalties.
Since the time horizon is set by the observation window and the  penalties govern the trade-off between control effort and agreement with the target laws, it is natural to ask whether the optimality system is well posed for arbitrary finite horizons and general penalties in this path-dependent setting.



For a more general class of path-dependent McKean--Vlasov control problems, \citet{cosso2023path} derive a dynamic programming principle and a path-dependent Bellman equation. 
For a squared MMD penalty on the joint law, Hafsi et al.~\cite[Theorem~5.15]{hafsi2026} and El Boustany et al.~\cite[Theorem~1]{elbou2026} establish solvability on sufficiently short horizons by a contraction argument. The admissible horizon depends on the penalty weights and the number of observation times. Both results assume that the spatial Jacobian of each column $\sigma^j$ of the diffusion matrix is a scalar multiple of the identity, $\partial_x\sigma^j(t,x)=\Sigma_j(t,x)I_d$.
Hafsi et al.~\cite[Theorem~5.7]{hafsi2026} also establish global solvability for a class of problems with linear dynamics and convex lifts of the cost functionals to $L^2$. Squared MMD is convex along mixtures of probability measures and can have a nonconvex lift to $L^2$, as illustrated by the Gaussian kernel \cite[Remark~3.5]{ma2026softbridge}.

\citet[Theorem~4.6]{ma2026softbridge} prove global solvability and uniqueness in law for a soft-constrained Schr\"odinger bridge with a penalty on the terminal marginal and a diffusion coefficient depending only on time.
They assume convexity of the terminal penalty along mixtures of laws. The control enters the drift through $\sigma(t)u$ and has running cost $\frac12|u|^2$. This relation between the control channel, the diffusion matrix and the running cost gives the relative entropy formulation and the conditional Gibbs representation. The corresponding Bellman equation at a fixed terminal law has the quadratic gradient term $-\frac12|\sigma(t)^\top\nabla v|^2$, which the Hopf--Cole transformation cancels against the corresponding term from the second-order operator.

Under our regularity assumptions, the two formulations are special cases of~\eqref{eq:cost}. The Schr\"odinger structure can persist with a state-dependent diffusion if the drift perturbation is $\sigma(t,\boldsymbol x)u$ and the control cost is $\frac12|u|^2$. The general drift and running cost in~\eqref{eq:cost} lead to Hamiltonians beyond this structure. For instance, Section~\ref{sec:affine} includes a nonquadratic running cost for which the minimized Hamiltonian has cubic growth in the adjoint variable. We establish the needed regularity directly for the nonlinear Bellman equations.

\citet{bai2026diffusion} study weak stochastic control on Wiener path space with convex law-dependent potentials. For a Brownian reference with deterministic initial state, their framework includes convex penalties on marginals and on the joint law of finitely many observations that satisfy their regularity assumptions.


Under uniform ellipticity and the regularity assumptions on the state-dependent diffusion coefficient and the flat derivatives of the penalties, together with conditions on the Hamiltonian including a growth bound on its spatial derivative, we prove that convexity along mixtures suffices for global well-posedness in the strong formulation. For every $T>0$ and every finite choice of nonnegative penalty weights, the control problem~\eqref{eq:cost} has a unique optimal control $u^*$, and the associated jump FBSDE~\eqref{eq:fbsde} has a unique square-integrable adapted solution $(X,Y,Z)$ with bounded $Y$. The constructed solution also has bounded $Z$.  \hypersetup{linkcolor=red}Under Assumption~\ref{ass:global}, uniqueness holds among all square-integrable solutions. The optimal control satisfies $u_t^*=\hat u(t,\boldsymbol X_t,Y_t)$, where $\hat u$ minimizes the Hamiltonian and $\boldsymbol X_t$ collects the recorded observations and the current state at time $t$. It is given by a feedback that is uniformly Lipschitz in $\boldsymbol X_t$ and grows at most 
linearly (Theorem~\ref{thm:main}). For squared MMD penalties with kernels satisfying Assumption~\ref{ass:mmd}, Corollary~\ref{cor:softpath} extends the local results of~\cite{hafsi2026,elbou2026} to arbitrary finite horizons and nonnegative penalty weights under our smoothness and ellipticity assumptions on a general state-dependent diffusion matrix.

The proof uses the potential structure of the associated mean field game \cite{lasrylions2007}.
We freeze the population path law and replace the two law penalties by their flat derivatives at that law to obtain the objective of a representative agent.
The agent's value function solves a sequence of Hamilton--Jacobi--Bellman (HJB) equations, with the past observations as parameters and matching conditions at the observation dates. Treating the recorded coordinates as parameters leaves a uniformly parabolic equation in dimension $d$ on each interval. 
The analytic step is to obtain a global Bernstein gradient bound and mixed derivative estimates in the current state and the frozen coordinates, uniformly in the population path law. These estimates yield a Lipschitz best-response feedback and its continuity with respect to the frozen law.
We obtain a  Nash  equilibrium by applying the Schauder--Tychonoff theorem to a compact set of path laws. This argument is related to the compactness method of Carmona and Delarue~\cite{cardel2013fbsde} for Markovian McKean--Vlasov FBSDEs.
In this potential game, convexity of the law penalties along
mixtures implies Lasry-Lions monotonicity of their
couplings (see Remark~\ref{monolasrylions}).
At the fixed point, the supporting inequalities from convexity along mixtures, combined with the HJB verification identity~\eqref{eq:verification}, identify the equilibrium control as the unique minimizer of the original objective.

At each intermediate observation date, the current state becomes a recorded coordinate and also remains the current-state argument of the next value function. Differentiating the matching condition~\eqref{eq:transmission} gives a derivative in each of these two arguments. Along the equilibrium trajectory, the derivative in the newly recorded coordinate is $Y_{s_j-}-Y_{s_j}$ and accounts for the effects of that observation on future dynamics and costs. Together with differentiation of the HJB equations between dates, this gives the jump FBSDE at the equilibrium law.
Conversely, every solution in the uniqueness class of Theorem~\ref{thm:main} induces a best response to its own path law and hence the same optimum.

We also establish global quadratic growth and gradient domination for the control cost $J$ when $b(t,x,u)=u$ and $\ell(t,x,u)=\frac12|u|^2$, 
with $\sigma$ depending on the current state and $F_t$ on its marginal law. Writing $Y^u$ for the adjoint associated with an admissible control $u$, we show that $u+Y^u$ represents the first variation of $J$ in $\mathbb H^2$ and that there exist constants $c,C>0$ such that
\[
 c\left(\|u-u^*\|_2^2+\E\sup_{0\le t\le T}|X_t^u-X_t^{u^*}|^2\right)
 \le J(u)-J(u^*)
 \le C\|u+Y^u\|_2^2.
\]
The left inequality is the quadratic growth, and the right one is the gradient domination. Together, the inequalities give an a posteriori bound on the cost gap and the squared control and state errors in terms of the Pontryagin residual $u+Y^u$ (Proposition~\ref{fv:prop:gradient}). For gradient domination, we fix the path law generated by an arbitrary admissible control. A linear backward stochastic differential equation (BSDE) with matching conditions bounds the discrepancy between its adjoint, augmented by the sensitivities to the recorded states, and the derivatives of the fixed-law value function in terms of the Pontryagin residual. In classical stochastic control, \citet{sethisika2024} use the control derivative of the Hamiltonian to study a gradient flow whose objective converges globally under convexity assumptions on the Hamiltonian and terminal cost.

The rest of the paper is organized as follows. Section~\ref{sec:model} formulates the control problem, states the assumptions and introduces the McKean--Vlasov FBSDE. Section~\ref{sec:proofwellposed} proves global well-posedness. Section~\ref{sec:examples} presents applications, including soft path-law constraints and drifts affine in the control. Section~\ref{sec:conclusion} concludes. Appendix~\ref{app:proofs} collects the auxiliary proofs, and Appendix~\ref{fv:sec:examples} provides further examples of smooth convex law penalties.

\section{Problem formulation}\label{sec:model}

\subsection{Notation}
Fix a dimension $d\ge1$ and a horizon $T>0$. Let $(\Omega,\F,\mathbb P)$ be a complete probability space supporting a $d$-dimensional Brownian motion $W$ and an $\R^d$-valued random variable $\xi$ with $\E|\xi|^2<\infty$, independent of $W$. We denote by $(\F_t)_{t\in[0,T]}$ the usual augmentation of the filtration generated by $\xi$ and $W$. We write $\Law(U)$ for the law of a random variable $U$ and $\delta_x$ for the unit point mass at $x$. We denote by $|\cdot|$ the Euclidean norm of a vector and by $\|\cdot\|$ the Frobenius norm of a matrix, and we define $A:B:=\operatorname{tr}(A^\top B)$.

For a normed space $E$, we denote by $\cP_2(E)$ the set of probability measures on $E$ with finite second moment, and we define the bounded Lipschitz metric
\[
 d_{\BL,E}(\mu,\nu):=
 \sup_{\|\varphi\|_\infty+\operatorname{Lip}(\varphi)\le1}
 \left|\int_E\varphi\,\dd(\mu-\nu)\right|,
\]
which metrizes weak convergence. We omit $E$ when the space is clear. We equip $C([0,T];\R^d)$ with the uniform norm. For $j\in\mathbb N_0:=\{0,1,2,\ldots\}$ and a Euclidean space $E$, $C_b^j(E)$ is the space of functions on $E$ whose derivatives up to order $j$ are bounded and continuous, with norm $\|\cdot\|_{C_b^j(E)}$ equal to the sum of their supremum norms; for $\alpha\in(0,1)$, $C_b^{j+\alpha}(E)$ and $\|\cdot\|_{C_b^{j+\alpha}(E)}$ add the $\alpha$-H\"older seminorm of the $j$-th derivatives. On product domains, $C_b^{k,j}$ and $C_b^{k,i,j}$ indicate the order of derivatives in corresponding variable.

For a Euclidean space $E$ of vectors or matrices, we denote by $\mathbb H^2(E)$ the space of $E$-valued $\F_t$-progressively measurable processes $V$ on $[0,T]$ with $$\norm V_2^2:=\E\int_0^T\|V_t\|^2\dt<\infty.$$ Let $\mathbb A\subset\R^q$ be a closed set of control values containing $0$. We denote by $\mathbb H^2(\mathbb A)$ the set of $\mathbb A$-valued processes in $\mathbb H^2(\R^q)$; controls are elements of $\mathbb H^2(\mathbb A)$.

\subsection{Law-penalized control problem and its optimality system}
Let $0<s_1<\cdots<s_m=T$, $m\in\mathbb N$, be the observation times, and define $s_0:=0$ and $\Pi:=(s_1,\ldots,s_m)$. We define the index $\delta^\Pi(t):=\max\{j<m:s_j\le t\}$ for $t\in[0,T]$. For a continuous process $X$, we denote by
\[
 \boldsymbol X_t:=(X_{s_1},\ldots,X_{s_{\delta^\Pi(t)}},X_t)
\]
the vector of the observations recorded up to $t$ and the current state. In particular, $\boldsymbol X_T=(X_{s_1},\ldots,X_{s_m})$.

For $u\in\mathbb H^2(\mathbb A)$, we denote by $X^u$ the solution of the state equation
\begin{equation}\label{eq:state}
 \dd X_t^u=b(t,\boldsymbol X_t^u,u_t)\dt+\sigma(t,\boldsymbol X_t^u)\dW_t,\qquad X_0^u=\xi,
\end{equation}
where $b(t,\cdot,\cdot):(\R^d)^{\delta^\Pi(t)+1}\times\mathbb A\to\R^d$ for $t\in[0,T]$ and $\sigma:[0,T]\times(\R^d)^{\delta^\Pi(t)+1}\to\R^{d\times d}$; under Assumption~\ref{ass:coefficients}\ref{ass:sigma}, \ref{ass:b0} and \ref{ass:drift} below, \eqref{eq:state} is well posed by Lemma~\ref{lem:state}.

Let $\ell(t,\cdot,\cdot):(\R^d)^{\delta^\Pi(t)+1}\times\mathbb A\to\R$ and $F_t:\cP_2((\R^d)^{\delta^\Pi(t)+1})\to\R$ for $t\in[0,T]$, and $G:\cP_2(\R^{md})\to\R$. We consider the control problem
\begin{equation}\label{eq:cost}
\inf_{u\in\mathbb H^2(\mathbb A)}J(u),\qquad
 J(u):=\E\int_0^T\ell(t,\boldsymbol X_t^u,u_t)\dt+\int_0^TF_t\bigl(\Law(\boldsymbol X_t^u)\bigr)\dt+G\bigl(\Law(\boldsymbol X_T^u)\bigr).
\end{equation}
In applications, $F_t(\mu)$ (resp.~$G(\eta)$) is typically a weighted discrepancy between $ \mu$ (resp.~$\eta$) and a target law $\rho$. Examples include the squared MMD in Section~\ref{sec:softpath} and smoothed versions of other discrepancies in Appendix~\ref{fv:sec:examples}.

We recall the flat derivative, also called the linear functional derivative; cf.~\cite[Definition~5.43]{carmonadelarue2018}. Let $E$ be a Euclidean space. A continuous function $\phi:\cP_2(E)\times E\to\R$ is a \emph{flat derivative} of a function $\Phi:\cP_2(E)\to\R$ if $\sup_{x\in E}|\phi(\mu,x)|/(1+|x|^2)$ is bounded on every subset of $\cP_2(E)$ with bounded second moments and, for all $\mu,\mu'\in\cP_2(E)$,
\begin{equation}\label{eq:flat}
 \Phi(\mu')-\Phi(\mu)=\int_0^1\int_E\phi\bigl((1-r)\mu+r\mu',x\bigr)(\mu'-\mu)(\dd x)\,\dd r.
\end{equation}

\begin{definition}\label{def:convex}
Let $E$ be a Euclidean space. A function $\Phi:\cP_2(E)\to\R$ is \emph{convex along mixtures of probability measures} if, for all $\mu,\mu'\in\cP_2(E)$ and $r\in[0,1]$,
\[
 \Phi\bigl((1-r)\mu+r\mu'\bigr)\le(1-r)\Phi(\mu)+r\Phi(\mu').
\]
\end{definition}

Definition~\ref{def:convex} differs from convexity in the sense of Lions, which is used in the probabilistic approach to McKean--Vlasov control~\cite{carmonadelarue2018} and in the global result~\cite[Theorem~5.7]{hafsi2026}: $\Phi$ is convex in the sense of Lions if, for all square-integrable $E$-valued random variables $U,U'$ on a common probability space and $r\in[0,1]$,
\[
 \Phi\bigl(\Law((1-r)U+rU')\bigr)\le(1-r)\Phi(\Law(U))+r\Phi(\Law(U')).
\]
This interpolates random variables, whereas Definition~\ref{def:convex} interpolates laws: $\Law((1-r)U+rU')$ depends on the coupling of $U$ and $U'$, while $(1-r)\Law(U)+r\Law(U')$ depends only on the marginals. For example, the squared MMD with a Gaussian kernel is convex along mixtures (Lemma~\ref{lem:mmd}), but not convex in the sense of Lions; see~\cite[Remark~3.5]{ma2026softbridge}.

\begin{remark}[Lasry--Lions monotonicity] \label{monolasrylions}
    When $\Phi$ admits a flat derivative $\phi$ as in~\eqref{eq:flat}, Definition~\ref{def:convex} is equivalent to the Lasry--Lions monotonicity condition
\begin{equation}\label{eq:ll-monotone}
 \int_E\bigl(\phi(\mu',x)-\phi(\mu,x)\bigr)(\mu'-\mu)(\dd x)\ge0,\qquad \mu,\mu'\in\cP_2(E),
\end{equation}
introduced in~\cite{lasrylions2007} to prove uniqueness of equilibria in mean field games. This equivalence can be found in \cite{graber2025Remarks, hofer2026Optimal}. We give a brief proof here. Fix $\mu,\mu'\in\cP_2(E)$ and define $\mu_r:=(1-r)\mu+r\mu'$ for $r\in[0,1]$, $\psi(r)$ $:=$ $\Phi(\mu_r)$. By the definition of the flat derivative, we have
\[
 \psi'(r)=\int_E\phi(\mu_r,x)(\mu'-\mu)(\dd x).
\]
By Definition \ref{def:convex}, $\Phi$ is convex along mixtures if and only if $\psi$ is convex on $ [0,1]$ for all $ \mu, \mu' \in \cP_2(E)$. For such each pair,  $\psi$ is convex if and only if $\psi'$ is nondecreasing, that is, for all $0\le r_1<r_2\le1$,
\[
 (r_2-r_1)\bigl(\psi'(r_2)-\psi'(r_1)\bigr)
 =\int_E\bigl(\phi(\mu_{r_2},x)-\phi(\mu_{r_1},x)\bigr)(\mu_{r_2}-\mu_{r_1})(\dd x)\ge0.
\]
Taking $r_1=0$ and $r_2=1$ yields~\eqref{eq:ll-monotone}; conversely,~\eqref{eq:ll-monotone} applied to $\mu_{r_1}$ and $\mu_{r_2}$ gives the inequality above and hence Definition~\ref{def:convex}.
\end{remark}

For $t\in[0,T]$, $\boldsymbol x\in(\R^d)^{\delta^\Pi(t)+1}$ and $y\in\R^d$, we define the Hamiltonian
\begin{equation}\label{eq:reducedhamiltonian}
 H(t,\boldsymbol x,y):=\inf_{u\in\mathbb A}\bigl\{b(t,\boldsymbol x,u)\cdot y+\ell(t,\boldsymbol x,u)\bigr\}.
\end{equation}
When the infimum in~\eqref{eq:reducedhamiltonian} is uniquely attained, denote the minimizer by $$\hat u(t,\boldsymbol x,y)\in \arg\min_{u\in \mathbb A} \left\{ b(t,\boldsymbol x,u) \cdot y + \ell(t,\boldsymbol x,u)\right\}.$$

Spatial derivatives of $\sigma$, $a$, $b$, $\ell$, $H$ and $f_t$ below are derivatives in $\boldsymbol x$. For $t\in[s_k,s_{k+1})$ and $1\le j\le k+1$, we denote by $\partial_{x_j}$ the derivative in the $j$-th component $x_j$ of $\boldsymbol x\in(\R^d)^{k+1}$, so that $\partial_{x_{k+1}}$ is the derivative in the current state. Since the dimension of $\boldsymbol x$ changes at the observation times, continuity in $t$ is required on each interval $[s_k,s_{k+1})$, $0\le k<m$.

We make the following assumptions in this paper.
\begin{assumption}\label{ass:coefficients}
There exist $\alpha\in(0,1)$, $L\ge0$, $\underline a>0$ and $C_0\ge0$ such that the following hold.
\begin{enumerate}[label=(\roman*)]
\item\label{ass:sigma} $\sigma:[0,T]\times(\R^d)^{\delta^\Pi(t)+1}\to\R^{d\times d}$ and its spatial derivatives up to order 4 are continuous, and for all $0\le k<m$, $s,t\in[s_k,s_{k+1})$ and $\boldsymbol x\in(\R^d)^{k+1}$,
\[
 \begin{gathered}
 \|\sigma(t,\cdot)\|_{C_b^{4+\alpha}}\le L,\qquad
 \|\sigma(t,\cdot)-\sigma(s,\cdot)\|_{C_b^1}\le L|t-s|^{1/2},\\
 a(t,\boldsymbol x):=\sigma(t,\boldsymbol x)\sigma(t,\boldsymbol x)^\top\ge \underline a I_d.
 \end{gathered}
\]
\item\label{ass:penalties} $F_t$ and $G$ are convex along mixtures of probability measures in the sense of Definition \ref{def:convex}, the map $t\mapsto F_t(\delta_0)$ is bounded and continuous, and $F_t$ and $G$ admit flat derivatives $f_t(\mu,\boldsymbol x)$ and $g(\eta,\boldsymbol x)$ with
\[
 \sup_{t,\mu}\|f_t(\mu,\cdot)\|_{C_b^{4+\alpha}}
 +\sup_\eta\|g(\eta,\cdot)\|_{C_b^{4+\alpha}}\le L
\]
and, for all $0\le k<m$, $s,t\in[s_k,s_{k+1})$, $\mu,\mu'\in\cP_2((\R^d)^{k+1})$ and $\eta,\eta'\in\cP_2(\R^{md})$,
\begin{equation}\label{eq:flatcont}
 \|f_t(\mu,\cdot)-f_s(\mu',\cdot)\|_{C_b^1}
 \le L\bigl(|t-s|^{1/2}+d_{\BL}(\mu,\mu')\bigr),\qquad
 \|g(\eta,\cdot)-g(\eta',\cdot)\|_{C_b^1}
 \le Ld_{\BL}(\eta,\eta').
\end{equation}
\item\label{ass:b0} $b(\cdot,\cdot,0)$ and its spatial derivatives up to order 4 are continuous, and for all $0\le k<m$ and $s,t\in[s_k,s_{k+1})$,
\[
 \begin{gathered}
 |b(t,0,0)|+\|\partial_{\boldsymbol x}b(t,\cdot,0)\|_{C_b^{3+\alpha}}\le L,\\
 |b(t,0,0)-b(s,0,0)|+\|\partial_{\boldsymbol x}b(t,\cdot,0)-\partial_{\boldsymbol x}b(s,\cdot,0)\|_{C_b^0}\le L|t-s|^{1/2}.
 \end{gathered}
\]
\item\label{ass:drift} $b$ and $\ell$ are measurable, continuous in $(\boldsymbol x,u)$ and differentiable in $\boldsymbol x$, with derivatives continuous in $(\boldsymbol x,u)$, and for all $t\in[0,T]$, $\boldsymbol x,\boldsymbol x'\in(\R^d)^{\delta^\Pi(t)+1}$ and $u\in\mathbb A$,
\[
 \begin{gathered}
 |b(t,\boldsymbol x,u)-b(t,\boldsymbol x,0)|\le C_0(1+|u|),\qquad
 |b(t,\boldsymbol x,u)-b(t,\boldsymbol x',u)|\le C_0(1+|u|)|\boldsymbol x-\boldsymbol x'|,\\
 -C_0\le\ell(t,\boldsymbol x,u)\le C_0(1+|\boldsymbol x|^2+|u|^2).
 \end{gathered}
\]
\item\label{ass:uhat} For every $(t,\boldsymbol x,y)$, the infimum in \eqref{eq:reducedhamiltonian} is attained at a unique point, and $\hat u$ is measurable and $\sup_{t\in[0,T]}|\hat u(t,0,0)|\le C_0$. For every $R>0$, $\hat u$ is Lipschitz continuous in $(\boldsymbol x,y)$ on $\{|y|\le R\}$, uniformly in $t$.
\item\label{ass:Hreg} The function $H$ and its derivatives in $(\boldsymbol x,y)$ up to order 4 exist and are continuous, and $H(\cdot,\cdot,0)$ is bounded. For every $R>0$, the derivatives of $H$ in $(\boldsymbol x,y)$ of order 2 to 4 are bounded on $\{(t,\boldsymbol{x},y):|y|\le R\}$.
\item\label{ass:bernstein} For all $(t,\boldsymbol x,y)$, $|\partial_{\boldsymbol x}H(t,\boldsymbol x,y)|\le C_0(1+|y|^2)$.
\end{enumerate}
\end{assumption}

In Assumption~\ref{ass:coefficients}\ref{ass:b0}, $b(t,\boldsymbol x,0)$ may grow linearly in $\boldsymbol x$, as for the Ornstein--Uhlenbeck drift $-\frac12x$. Condition~\ref{ass:drift} allows the spatial Lipschitz constant of $b$ to grow with the control, as in the case $b(t,\boldsymbol x,u)=b(t,\boldsymbol x,0)+\sigma(t,\boldsymbol x)u$. Condition~\ref{ass:uhat} is a condition on the minimizer itself rather than on $H$. It holds, for instance, if $\mathbb A=\R^q$, $b$ is affine in $u$, and $\ell$ is uniformly convex in $u$; see Corollary~\ref{cor:affine} for the precise conditions. Condition~\ref{ass:bernstein} yields the gradient bound in Lemma~\ref{lem:hjb}. It also makes $\partial_{\boldsymbol x}H(t,\boldsymbol X_t,Y_t)$ integrable when $Y$ is square integrable, so that the jump condition in~\eqref{eq:fbsde} below is well defined for square-integrable solutions. Since the spatial derivatives of $\sigma$ are bounded, the term $\partial_{x_j}\sigma(t,\boldsymbol X_t):Z_t$ in this jump condition is integrable when $Z\in\mathbb H^2(\R^{d\times d})$. The following lemma gives the well-posedness of the state $ X$.

\begin{lemma}\label{lem:state}
Suppose Assumption~\ref{ass:coefficients}\ref{ass:sigma}, \ref{ass:b0} and \ref{ass:drift} hold. Then, for every $u\in\mathbb H^2(\mathbb A)$, the equation \eqref{eq:state} has a unique strong solution, and $\E\sup_{t\le T}|X_t^u|^2<\infty$.
\end{lemma}

This lemma is proved in Appendix~\ref{app:state} by applying the well-posedness result of \cite{gyongyKrylov1980} successively on the intervals between observation times.

Under Assumption~\ref{ass:coefficients}, $J(u)$ is finite for every $u\in\mathbb H^2(\mathbb A)$: the running cost $\ell(t,\boldsymbol X_t^u,u_t)$ is integrable by~\ref{ass:drift} and Lemma~\ref{lem:state}, while~\ref{ass:penalties} and~\eqref{eq:flat} show that $G$ is bounded and that $t\mapsto F_t(\Law(\boldsymbol X_t^u))$ is bounded and continuous on each $[s_k,s_{k+1})$. Fix $t\in[0,T)$ and $\mu,\mu'\in\cP_2((\R^d)^{\delta^\Pi(t)+1})$. By Definition~\ref{def:convex} and~\eqref{eq:flat}, for $\eps\in(0,1]$,
\[
 \begin{aligned}
 F_t(\mu')-F_t(\mu)
 &\ge\frac{F_t((1-\eps)\mu+\eps\mu')-F_t(\mu)}\eps\\
 &=\int_0^1\int_{(\R^d)^{\delta^\Pi(t)+1}} f_t\bigl((1-r\eps)\mu+r\eps\mu',\boldsymbol x\bigr)(\mu'-\mu)(\dd \boldsymbol x)\,\dd r.
 \end{aligned}
\]
Since $d_{\BL}((1-r\eps)\mu+r\eps\mu',\mu)=r\eps\,d_{\BL}(\mu',\mu)\le2r\eps$, it follows from~\eqref{eq:flatcont} that $f_t((1-r\eps)\mu+r\eps\mu',\cdot)$ converges uniformly to $f_t(\mu,\cdot)$ as $\eps\downarrow0$. Letting $\eps\downarrow0$ and arguing similarly for $G$, we obtain, for all $t\in[0,T)$, $\mu,\mu'\in\cP_2((\R^d)^{\delta^\Pi(t)+1})$ and $\eta,\eta'\in\cP_2(\R^{md})$,
\begin{equation}\label{eq:convex}
\begin{aligned}
 F_t(\mu')-F_t(\mu)&\ge\int_{(\R^d)^{\delta^\Pi(t)+1}} f_t(\mu,\boldsymbol x)(\mu'-\mu)(\dd \boldsymbol x),\\
 G(\eta')-G(\eta)&\ge\int_{\R^{md}} g(\eta,\boldsymbol x)(\eta'-\eta)(\dd\boldsymbol x).
\end{aligned}
\end{equation}

By Assumption~\ref{ass:coefficients}\ref{ass:uhat}, for each $R>0$ there exists $C_R$ such that, for $|y|\le R$,
\begin{equation}\label{eq:uhatbounds}
|\hat u(t,\boldsymbol x,y)|\le C_R(1+|\boldsymbol x|).
\end{equation}
For fixed $t$, the minimizer thus stays in a compact set locally in $(\boldsymbol x,y)$. By the continuity and Assumption~\ref{ass:coefficients}\ref{ass:drift}, we have 
\begin{equation}\label{eq:envelope}
 \partial_yH(t,\boldsymbol x,y)=b\bigl(t,\boldsymbol x,\hat u(t,\boldsymbol x,y)\bigr),\quad
 \partial_{\boldsymbol x}H(t,\boldsymbol x,y)=\partial_{\boldsymbol x}b\bigl(t,\boldsymbol x,\hat u(t,\boldsymbol x,y)\bigr)^\top y+\partial_{\boldsymbol x}\ell\bigl(t,\boldsymbol x,\hat u(t,\boldsymbol x,y)\bigr).
\end{equation}
As an infimum of affine functions of $y$, $H(t,\boldsymbol x,\cdot)$ is concave, so by Assumption \ref{ass:coefficients}\ref{ass:Hreg},
\begin{equation}\label{eq:concavity}
 H(t,\boldsymbol x,y)-\partial_yH(t,\boldsymbol x,y)\cdot y\ge H(t,\boldsymbol x,0)\ge-\sup_{t,\boldsymbol x}|H(t,\boldsymbol x,0)|>-\infty.
\end{equation}
We also introduce the following stronger assumption.
\begin{assumption}\label{ass:global}
For all $t\in[0,T]$, $\boldsymbol x\in(\R^d)^{\delta^\Pi(t)+1}$ and $y\in\R^d$, with $C_0$ as in Assumption~\ref{ass:coefficients},
\begin{equation}\label{eq:globalbounds}
 \|\partial_{yy}H(t,\boldsymbol x,y)\|+\|\partial_{\boldsymbol xy}H(t,\boldsymbol x,y)\|\le C_0.
\end{equation}
\end{assumption}

Assumption~\ref{ass:global} extends the bounds on $\partial_{yy}H$ and $\partial_{\boldsymbol xy}H$ in Assumption~\ref{ass:coefficients}\ref{ass:Hreg} to all $y$. It is used only for uniqueness among square-integrable solutions, in Lemma~\ref{lem:uniqueness} and Theorem~\ref{thm:main}\ref{thm:main-c}.

Under Assumption~\ref{ass:coefficients}\ref{ass:penalties}, $\partial_{\boldsymbol x}f_t$ and $\partial_{\boldsymbol x}g$ are bounded and jointly continuous. By~\cite[Proposition~5.48]{carmonadelarue2018}, $\partial_{\boldsymbol x}f_t(\mu,\cdot)$ is the Lions derivative of $F_t$ at $\mu$, and $\partial_{\boldsymbol x}g(\eta,\cdot)=(\partial_{x_1}g(\eta,\cdot),\ldots,\partial_{x_m}g(\eta,\cdot))$ is the Lions derivative of $G$ at $\eta$. For $z\in\R^{d\times d}$, $t\in[s_k,s_{k+1})$ and $1\le j\le k+1$, we denote by $\partial_{x_j}\sigma(t,\boldsymbol x):z$ the gradient in $x_j$ of $\boldsymbol x\mapsto\sigma(t,\boldsymbol x):z$.

We consider the McKean--Vlasov FBSDE
\begin{equation}\label{eq:fbsde}
 \begin{array}{rcl}
 \dd X_t & =&   b(t,\boldsymbol X_t,\hat u_t)\dt+\sigma(t,\boldsymbol X_t)\dW_t,
        \qquad X_0=\xi,\\
 \dd Y_t &= & -\bigl[\partial_{x_{k+1}}H(t,\boldsymbol X_t,Y_t)+\partial_{x_{k+1}}\sigma(t,\boldsymbol X_t):Z_t\\
  &  &\qquad +\partial_{x_{k+1}}f_t(\Law(\boldsymbol X_t),\boldsymbol X_t)\bigr]\dt+Z_t\dW_t,
       \quad t\in(s_k,s_{k+1}),\ 0\le k<m,\\
 Y_{s_j}-Y_{s_j-}& =& -\E\Bigl[\partial_{x_j}g(\Law(\boldsymbol X_T),\boldsymbol X_T)   + \int_{s_j}^T\bigl(\partial_{x_j}H(t,\boldsymbol X_t,Y_t) \\
& & \quad \quad  +\partial_{x_j}\sigma(t,\boldsymbol X_t):Z_t 
+\partial_{x_j}f_t(\Law(\boldsymbol X_t),\boldsymbol X_t)\bigr)\dt\Bigm|\F_{s_j}\Bigr], \quad 1\le j<m,\\
 Y_T & = &\partial_{x_m}g(\Law(\boldsymbol X_T),\boldsymbol X_T),
 \end{array}
\end{equation}
where $\hat u_t:=\hat u(t,\boldsymbol X_t,Y_t)$ and $Y_{s_j-}$ denotes the left limit at $s_j$. Since $\sigma$ is uncontrolled, by~\eqref{eq:envelope}, \eqref{eq:fbsde} is the optimality system of the stochastic maximum principle for~\eqref{eq:cost}; see~\cite{carmonadelarue2015,carmonadelarue2018} for marginal-law penalties and~\cite{hafsi2026} for discrete path dependence. The sufficient condition~\cite[Theorem~4.7]{hafsi2026} requires the terminal cost and the Hamiltonian of the stochastic maximum principle to be convex in the sense of~\cite[Definition~3.2]{hafsi2026}, jointly in the state, the law and the control; this fails in general, for instance through $\boldsymbol x\mapsto\sigma(t,\boldsymbol x):Z_t$.

We call an $\F_t$-adapted solution $(X,Y,Z)$ of \eqref{eq:fbsde} \emph{square integrable} if $X$ is continuous, $Y$ is c\`adl\`ag with $Y_{T-}=Y_T$, $\E\sup_{t\le T}(|X_t|^2+|Y_t|^2)<\infty$ and $Z\in\mathbb H^2(\R^{d\times d})$. Uniqueness of solutions is understood up to indistinguishability for $(X,Y)$ and in $\mathbb H^2(\R^{d\times d})$ for $Z$, and uniqueness of controls up to $\dd t\otimes\dd\mathbb P$-null sets.

\section{Global well-posedness}\label{sec:proofwellposed}

We solve the mean field control problem~\eqref{eq:cost} and the McKean--Vlasov FBSDE~\eqref{eq:fbsde} through an associated mean field game. We first fix the path law $\mathbb Q$ of the population and solve the HJB equations of a single agent facing $\mathbb Q$ (Lemmas~\ref{lem:hjb} and~\ref{lem:frozenproblem}). The optimal feedback of the agent defines the state equation~\eqref{eq:frozenstate}, whose solution has path law $\mathcal T\mathbb Q$ (see~\eqref{eq:bestresponse}). We then find a fixed point $\mathbb Q^*=\mathcal T\mathbb Q^*$, which is an equilibrium of the game (Lemma~\ref{lem:fixedpoint}), and convexity along mixtures shows that this equilibrium solves the mean field control problem (Lemma~\ref{lem:optimal}). Finally, Lemmas~\ref{lem:construction} and~\ref{lem:uniqueness} establish existence and uniqueness for~\eqref{eq:fbsde}, yielding Theorem~\ref{thm:main}.

For $\mathbb Q\in\cP_2(C([0,T];\R^d))$, we denote by $\mathbb Q_t$, $\mathbb Q_{\Pi(t)}$ and $\mathbb Q_\Pi$ the laws of $\omega_t$, $\boldsymbol\omega_t:=(\omega_{s_1},\ldots,\omega_{s_{\delta^\Pi(t)}},\omega_t)$ and $(\omega_{s_1},\ldots,\omega_{s_m})$ under $\mathbb Q$, where $\Pi(t):=(s_1,\ldots,s_{\delta^\Pi(t)},t)$, so that $\mathbb Q_{\Pi(T)}=\mathbb Q_\Pi$. Given the path law $\mathbb Q$ of a population, we linearize the penalties at $\mathbb Q$: we replace $F_t$ and $G$ by their flat derivatives $f_t(\mathbb Q_{\Pi(t)},\cdot)$ and $g(\mathbb Q_\Pi,\cdot)$. Then a single agent minimizes
\begin{equation}\label{eq:frozencost}
 J_{\mathbb Q}(u):=\E\left[\int_0^T\bigl(\ell(t,\boldsymbol X_t^u,u_t)+f_t(\mathbb Q_{\Pi(t)},\boldsymbol X_t^u)\bigr)\dt+g(\mathbb Q_\Pi,\boldsymbol X_T^u)\right],\qquad u\in\mathbb H^2(\mathbb A).
\end{equation}
This defines an associated \emph{mean field game}: given $\mathbb Q$, the agent solves $\inf_{u\in\mathbb H^2(\mathbb A)}J_{\mathbb Q}(u)$, and a path law $\mathbb Q^*\in\cP_2(C([0,T];\R^d))$ is a  \emph{Nash equilibrium} if there exists a minimizer $u^*\in\mathbb H^2(\mathbb A)$ of $J_{\mathbb Q^*}$ such that $\Law(X^{u^*})=\mathbb Q^*$.

Applying \eqref{eq:ll-monotone} to the joint laws entering \(f_t\) and \(g\), we obtain Lasry–Lions monotonicity of the coupling in \eqref{eq:frozencost} on population path laws by integrating the running pairing over time and adding the terminal pairing.

In general, the equilibrium of a mean field game and the optimum of a mean field control problem differ: the agents of the game ignore their influence on the population, while the planner of the control problem accounts for it~\cite{carmona2013}. Here the agents pay the marginal costs $f_t$ and $g$, and so internalize this influence. The game is then potential in the sense of Lasry and Lions~\cite{lasrylions2007}: its equations are the first-order optimality conditions of \eqref{eq:cost}. Convexity along mixtures implies that any Nash equilibrium control is a global minimizer of \eqref{eq:cost}. 
If $\mathbb Q^*$ is an equilibrium and $u^*$ is the optimal control of the agent at $\mathbb Q^*$, then~\eqref{eq:convex} gives, for every $u\in\mathbb H^2(\mathbb A)$,
\[
 J(u)-J(u^*)\ge J_{\mathbb Q^*}(u)-J_{\mathbb Q^*}(u^*)\ge0,
\]
since, up to a constant, $J-J_{\mathbb Q^*}$ is the sum of the remainders of the linearizations of the penalties at $\mathbb Q^*$, which are nonnegative and vanish at $u^*$. Since $J_{\mathbb Q^*}$ has a unique minimizer, the equilibrium control is the unique optimal control. 
Conversely, each solution of~\eqref{eq:fbsde} in the uniqueness class of Theorem~\ref{thm:main} follows the best response to its own path law. Lemma~\ref{lem:uniqueness} shows that it coincides with the solution associated with $u^*$. Neither step uses a contraction inequality: the equilibrium comes from compactness and its optimality from convexity, which is why the horizon in Theorem~\ref{thm:main} is arbitrary.

Throughout the proofs, $C$ denotes a generic constant, depending only on the quantities in the statement being proved, whose value may change from line to line.

\subsection{The HJB equations of the agent}
Minimizing $J_{\mathbb Q}$ \eqref{eq:frozencost} is a classical stochastic control problem with running cost $\ell(t,\boldsymbol x,u)+f_t(\mathbb Q_{\Pi(t)},\boldsymbol x)$ and terminal cost $g(\mathbb Q_\Pi,\boldsymbol x)$. Since the coefficients and the terminal cost depend on the observations $X_{s_1},\ldots,X_{s_m}$, the value function at time $t\in[s_k,s_{k+1})$, $0\le k<m$, depends on the current state $x\in\R^d$ and on the observations $\boldsymbol h_k=(h_1,\ldots,h_k)\in\R^{kd}$ already recorded:
\[
 \begin{aligned}
 v_k^{\mathbb Q}(t,\boldsymbol h_k,x)
 =\inf_{u\in\mathbb H^2(\mathbb A)}\E\Bigl[&\int_t^T\bigl(\ell(s,\boldsymbol X_s^{t,x,u},u_s)+f_s(\mathbb Q_{\Pi(s)},\boldsymbol X_s^{t,x,u})\bigr)\ds\\
   &+g\bigl(\mathbb Q_\Pi,\boldsymbol h_k,X_{s_{k+1}}^{t,x,u},\ldots,X_{s_m}^{t,x,u}\bigr)\Bigr],
 \end{aligned}
\]
where $X^{t,x,u}$ solves $\dd X_s=b(s,\boldsymbol X_s,u_s)\ds+\sigma(s,\boldsymbol X_s)\dW_s$ on $[t,T]$ with $X_t=x$, and $\boldsymbol X_s:=(\boldsymbol h_k,X_{s_{k+1}},\ldots,X_{s_{\delta^\Pi(s)}},X_s)$. Between observation times, the recorded observations do not move, so $v_k^{\mathbb Q}$ is expected to solve the HJB equation
\[
 \begin{aligned}
 &\partial_tv+\tfrac12a(t,(\boldsymbol h_k,x)):\partial_{xx}v\\
 &\quad{}+\inf_{u\in\mathbb A}\bigl\{b(t,(\boldsymbol h_k,x),u)\cdot\partial_xv+\ell(t,(\boldsymbol h_k,x),u)\bigr\}+f_t(\mathbb Q_{\Pi(t)},(\boldsymbol h_k,x))=0
 \end{aligned}
\]
for $ s_k\le t<s_{k+1}$, in which $\boldsymbol h_k$ enters only as a parameter. By \eqref{eq:reducedhamiltonian} and Assumption~\ref{ass:coefficients}\ref{ass:uhat}, the infimum equals $H(t,(\boldsymbol h_k,x),\partial_xv)$ and is uniquely attained at $u=\hat u(t,(\boldsymbol h_k,x),\partial_xv)$, which gives \eqref{eq:hjb} below and, for $t\in[s_k,s_{k+1})$, the optimal feedback
\[
 u_t=\hat u\bigl(t,(\boldsymbol h_k,x),\partial_xv_k^{\mathbb Q}(t,\boldsymbol h_k,x)\bigr)\quad\text{at}\quad\boldsymbol h_k=(X_{s_1},\ldots,X_{s_k}),\ x=X_t.
\]
Since $\boldsymbol X_t=(X_{s_1},\ldots,X_{s_k},X_t)$ for these $t$, we write this as $u_t=\hat u\bigl(t,\boldsymbol X_t,\partial_{x_{k+1}}v_k^{\mathbb Q}(t,\boldsymbol X_t)\bigr)$, where $\partial_{x_{k+1}}$ is the derivative in the last of the $k+1$ variables of $v_k^{\mathbb Q}$. We first study the equation on a single interval, with a parameter $h\in\R^p$ standing for the recorded observations.

\begin{lemma}\label{lem:hjb}
Suppose Assumption~\ref{ass:coefficients} holds. Let $0\le k<m$, $[\underline t,\bar t]\subset[s_k,s_{k+1}]$, $p:=kd$, $U:[\underline t,\bar t]\times\R^p\times\R^d\to\R$ and $\Phi:\R^p\times\R^d\to\R$. Assume that $U\in C_b^{0,4}([\underline t,\bar t]\times\R^{p+d})$
and $\Phi\in C_b^4(\R^{p+d})$. Then the equation
\begin{equation}\label{eq:paramhjb}
 \begin{gathered}
 \partial_tv+\tfrac12a(t,(h,x)):\partial_{xx}v+H(t,(h,x),\partial_xv)+U(t,h,x)=0\quad\text{on }[\underline t,\bar t)\times\R^p\times\R^d,\\
 v(\bar t,h,x)=\Phi(h,x),
 \end{gathered}
\end{equation}
admits a classical solution
$v\in C_b^{0,4}([\underline t,\bar t]\times\R^{p+d})$. This solution is unique among bounded
classical solutions with bounded $\partial_xv$.
\end{lemma}
The proof, given in Appendix~\ref{app:hjb}, uses a Bernstein gradient bound and estimates of the mixed derivatives in $(h,x)$.

Let $\mathbb Q\in\cP_2(C([0,T];\R^d))$. For $0\le k\le m-1$, we construct backwards functions $v_k^{\mathbb Q}(t,\boldsymbol h_k,x)$, $t\in[s_k,s_{k+1})$, where $\boldsymbol h_k=(h_1,\ldots,h_k)\in\R^{kd}$ stands for the observations recorded at $s_1,\ldots,s_k$ and $x\in\R^d$ for the current state; for $k=0$, there is no recorded observation and $v_0^{\mathbb Q}=v_0^{\mathbb Q}(t,x)$. For every $\boldsymbol h_k\in\R^{kd}$, $v_k^{\mathbb Q}$ solves
\begin{equation}\label{eq:hjb}
 \begin{aligned}
 &\partial_tv_k^{\mathbb Q}+\tfrac12a(t,(\boldsymbol h_k,x)):\partial_{xx}v_k^{\mathbb Q}+H(t,(\boldsymbol h_k,x),\partial_xv_k^{\mathbb Q})\\
 &\quad{}+f_t(\mathbb Q_{\Pi(t)},(\boldsymbol h_k,x))=0,\qquad s_k\le t<s_{k+1},
 \end{aligned}
\end{equation}
with
\begin{align}
 v_{m-1}^{\mathbb Q}(T,\boldsymbol h_{m-1},x)&=g(\mathbb Q_\Pi,\boldsymbol h_{m-1},x),\label{eq:terminal}\\
 v_{k-1}^{\mathbb Q}(s_k-,\boldsymbol h_{k-1},x)&=v_k^{\mathbb Q}\bigl(s_k,(\boldsymbol h_{k-1},x),x\bigr),\qquad 1\le k<m.
 \label{eq:transmission}
\end{align}
In \eqref{eq:terminal}, we identify $(\boldsymbol h_{m-1},x)$ with $\boldsymbol x\in\R^{md}$. In \eqref{eq:transmission}, the current state $x$ at time $s_k$ becomes the $k$-th recorded observation.

\begin{lemma}\label{lem:frozenproblem}
Suppose Assumption~\ref{ass:coefficients} holds, and let $\mathbb Q\in\cP_2(C([0,T];\R^d))$. Then \eqref{eq:hjb}--\eqref{eq:transmission} admit solutions
$v_k^{\mathbb Q}$ satisfying
$\max_{0\le k<m}\sup_{t\in[s_k,s_{k+1})}
\|v_k^{\mathbb Q}(t,\cdot)\|_{C_b^4(\R^{(k+1)d})}\le C$
for some constant $C$ independent of $\mathbb Q$. In particular, there exists $\bar C\ge0$, independent of $\mathbb Q$, such that
\begin{equation}\label{eq:feedbackbounds}
 \sup_{k,t,\boldsymbol h_k,x}
 \bigl(|\partial_xv_k^{\mathbb Q}(t,\boldsymbol h_k,x)|+\|\partial_{(\boldsymbol h_k,x)}\partial_xv_k^{\mathbb Q}(t,\boldsymbol h_k,x)\|\bigr)\le\bar C.
\end{equation}
The equation
\begin{equation}\label{eq:frozenstate}
 \begin{gathered}
 \dd X_t^{\mathbb Q}=b(t,\boldsymbol X_t^{\mathbb Q},u_t^{\mathbb Q})\dt+\sigma(t,\boldsymbol X_t^{\mathbb Q})\dW_t,\qquad X_0^{\mathbb Q}=\xi,\\
 u_t^{\mathbb Q}:=\hat u\bigl(t,\boldsymbol X_t^{\mathbb Q},\partial_{x_{k+1}}v_k^{\mathbb Q}(t,\boldsymbol X_t^{\mathbb Q})\bigr),
 \end{gathered}
\end{equation}
where $k$ is such that $t\in[s_k,s_{k+1})$, has a unique strong solution, and $u^{\mathbb Q}$ is the unique minimizer of $J_{\mathbb Q}$ in $\mathbb H^2(\mathbb A)$. Finally, if $(\mathbb Q^n)_{n\in\mathbb N}\subset\cP_2(C([0,T];\R^d))$ is such that
\begin{equation}\label{eq:frozenhyp}
 \lim_{n\to\infty}\Bigl(\sup_{t<T,\boldsymbol x}|f_t(\mathbb Q_{\Pi(t)}^n,\boldsymbol x)-f_t(\mathbb Q_{\Pi(t)},\boldsymbol x)|
 +\sup_{\boldsymbol x}|g(\mathbb Q_\Pi^n,\boldsymbol x)-g(\mathbb Q_\Pi,\boldsymbol x)|\Bigr)=0,
\end{equation}
then, for each $0\le k<m$,
\begin{equation}\label{eq:frozengrad}
 \lim_{n\to\infty}\sup_{t\in[s_k,s_{k+1}),\,\boldsymbol h_k,\,x}
 |\partial_xv_k^{\mathbb Q^n}(t,\boldsymbol h_k,x)-\partial_xv_k^{\mathbb Q}(t,\boldsymbol h_k,x)|=0.
\end{equation}
\end{lemma}
The proof, given in  Appendix~\ref{app:frozen}, applies Lemma~\ref{lem:hjb} backward over the intervals $[s_k,s_{k+1}]$. The bounds in~\eqref{eq:feedbackbounds}, which are uniform in $\mathbb Q$, together with~\eqref{eq:frozengrad} make the best response map defined below continuous.

\subsection{Fixed point and optimal control}
For $\mathbb Q\in\cP_2(C([0,T];\R^d))$, let $X^{\mathbb Q}$ be the optimal state of the agent facing the path law $\mathbb Q$, that is, the solution of \eqref{eq:frozenstate} in Lemma~\ref{lem:frozenproblem}. We define the best response map $\mathcal T:\cP_2(C([0,T];\R^d))\to\cP_2(C([0,T];\R^d))$ by
\begin{equation}\label{eq:bestresponse}
 \mathcal T\mathbb Q:=\Law(X^{\mathbb Q}).
\end{equation}
A fixed point of $\mathcal T$ is an equilibrium of the mean field game associated with~\eqref{eq:frozencost}.

\begin{lemma}\label{lem:fixedpoint}
Suppose Assumption~\ref{ass:coefficients} holds. Then $\mathcal T$ has a fixed point $\mathbb Q^*\in\cP_2(C([0,T];\R^d))$.
\end{lemma}

This lemma is proved in Appendix~\ref{app:fixedpoint} by applying the Schauder--Tychonoff theorem to $\mathcal T$ on a convex compact set of path laws.

In the rest of this section, $\mathbb Q^*$ denotes a fixed point of $\mathcal T$ given by Lemma~\ref{lem:fixedpoint}. We define
\begin{equation}\label{eq:optimal control}
u_t^*:=\hat u\bigl(t,\boldsymbol X_t^{\mathbb Q^*},\partial_{x_{k+1}}v_k^{\mathbb Q^*}(t,\boldsymbol X_t^{\mathbb Q^*})\bigr)
\end{equation}
 for $t\in[s_k,s_{k+1})$ and $0\le k<m$, so $u^*=u^{\mathbb Q^*}$ (see \eqref{eq:frozenstate}). By pathwise uniqueness for~\eqref{eq:state} (Lemma~\ref{lem:state}) and since $\mathbb Q^*$ is a fixed point of the map $\mathcal T$ in~\eqref{eq:bestresponse}, we have
\[
 X^{u^*}=X^{\mathbb Q^*},\qquad
 \Law(\boldsymbol X_t^{u^*})=\mathbb Q_{\Pi(t)}^*,\qquad \Law(\boldsymbol X_T^{u^*})=\mathbb Q_\Pi^*,
\]
and $u^*$ minimizes $J_{\mathbb Q^*}$ by Lemma~\ref{lem:frozenproblem}.

For $\mathbb Q\in\cP_2(C([0,T];\R^d))$, we define
\[
 c_{\mathbb Q}:=\int_0^T\left[F_t(\mathbb Q_{\Pi(t)})-\int_{(\R^d)^{\delta^\Pi(t)+1}} f_t(\mathbb Q_{\Pi(t)},\boldsymbol x)\,\mathbb Q_{\Pi(t)}(\dd \boldsymbol x)\right]\dt
       +G(\mathbb Q_\Pi)-\int_{\R^{md}} g(\mathbb Q_\Pi,\boldsymbol x)\,\mathbb Q_\Pi(\dd\boldsymbol x),
\]
and, for $t\in[0,T]$, $\mu,\nu\in\cP_2((\R^d)^{\delta^\Pi(t)+1})$ and $\eta,\zeta\in\cP_2(\R^{md})$,
\begin{align*}
 B_{F_t}(\mu,\nu)
 &:=F_t(\mu)-F_t(\nu)-\int_{(\R^d)^{\delta^\Pi(t)+1}} f_t(\nu,\boldsymbol x)(\mu-\nu)(\dd \boldsymbol x),\\
 B_G(\eta,\zeta)
 &:=G(\eta)-G(\zeta)-\int_{\R^{md}} g(\zeta,\boldsymbol x)(\eta-\zeta)(\dd\boldsymbol x).
\end{align*}
Substituting~\eqref{eq:cost} and~\eqref{eq:frozencost} gives
\begin{equation}\label{eq:linearization}
J(u)=J_{\mathbb Q}(u)+c_{\mathbb Q}+\int_0^T B_{F_t}\bigl(\Law(\boldsymbol X_t^u),\mathbb Q_{\Pi(t)}\bigr)\dt
                            +B_G\bigl(\Law(\boldsymbol X_T^u),\mathbb Q_\Pi\bigr),
\end{equation}
where every term is finite by Assumption~\ref{ass:coefficients} \ref{ass:penalties}, \ref{ass:drift} and Lemma~\ref{lem:state}.

We obtain, for every $u\in\mathbb H^2(\mathbb A)$,
\[
 J(u)\ge J_{\mathbb Q^*}(u)+c_{\mathbb Q^*}
        \ge J_{\mathbb Q^*}(u^*)+c_{\mathbb Q^*}=J(u^*),
\]
where the first inequality holds because $B_{F_t}$ and $B_G$ are nonnegative by the convexity of $F_t$ and $G$ (see~\eqref{eq:convex}), the second holds because $u^*$ minimizes $J_{\mathbb Q^*}$, and the last equality holds because the last two terms of~\eqref{eq:linearization} vanish at $\mathbb Q=\mathbb Q^*$ and $u=u^*$ since $\mathbb Q^*$ is a fixed point of $\mathcal T$ in~\eqref{eq:bestresponse}. Moreover, if $\bar u\in\mathbb H^2(\mathbb A)$ is also optimal for~\eqref{eq:cost}, then both inequalities above are equalities for $u=\bar u$, so $\bar u$ minimizes $J_{\mathbb Q^*}$, and Lemma~\ref{lem:frozenproblem} gives $\bar u=u^*$ up to $\dd t\otimes\dd\mathbb P$-null sets. This proves the following lemma.

\begin{lemma}\label{lem:optimal}
Suppose Assumption~\ref{ass:coefficients} holds. Then $u^*$ is the unique optimal control of~\eqref{eq:cost} in $\mathbb H^2(\mathbb A)$.
\end{lemma}

\subsection{Construction of the FBSDE solution}
\begin{lemma}\label{lem:construction}
Suppose Assumption~\ref{ass:coefficients} holds. We define $X:=X^{u^*}$ and, for $t\in[s_k,s_{k+1})$ and $0\le k<m$,
\[
 Y_t:=\partial_{x_{k+1}}v_k^{\mathbb Q^*}(t,\boldsymbol X_t),\qquad
 Z_t:=\partial_{x_{k+1}}^2v_k^{\mathbb Q^*}(t,\boldsymbol X_t)\sigma(t,\boldsymbol X_t),
\]
and $Y_T:=\partial_{x_m}g(\mathbb Q_\Pi^*,\boldsymbol X_T)$.\footnote{A predictable version of $Z$ is obtained by taking left limits at the observation times; this does not affect its stochastic integral.}
Then $(X,Y,Z)$ is a solution of \eqref{eq:fbsde} with $u_t^*=\hat u(t,\boldsymbol X_t,Y_t)$, $\E\sup_{t\le T}|X_t|^2<\infty$, and $|Y_t|,\|Z_t\|$ being bounded for $\dd t\otimes\dd\mathbb P$-almost every $(t,\omega)$. Moreover, the function $y^*(t,\boldsymbol x):=\partial_{x_{k+1}}v_k^{\mathbb Q^*}(t,\boldsymbol x)$, $t\in[s_k,s_{k+1})$, $\boldsymbol x\in(\R^d)^{k+1}$, $0\le k<m$, is bounded and Lipschitz continuous in $\boldsymbol x$, uniformly in $t$. The feedback $\hat u(t,\boldsymbol x,y^*(t,\boldsymbol x))$ is Lipschitz and has linear growth in $\boldsymbol x$, uniformly in $ t$.
\end{lemma}
\begin{proof}
By the definition \eqref{eq:optimal control}, we have $u_t^*=\hat u(t,\boldsymbol X_t,Y_t)$. 
On each interval $(s_k,s_{k+1})$, we apply It\^o's formula to $Y_t=\partial_{x_{k+1}}v_k^{\mathbb Q^*}(t,\boldsymbol X_t)$. Note that the drift of $X$ is $\partial_yH(t,\boldsymbol X_t,Y_t)$ by~\eqref{eq:envelope}.
\[
 \begin{aligned}
 \dd Y_t
 &=\Bigl(\partial_t\partial_{x_{k+1}}v_k^{\mathbb Q^*}(t,\boldsymbol X_t)+\partial_{x_{k+1}}^2v_k^{\mathbb Q^*}(t,\boldsymbol X_t)\,\partial_yH(t,\boldsymbol X_t,Y_t)\\
 &\qquad+\tfrac12a(t,\boldsymbol X_t):\partial_{x_{k+1}}^3v_k^{\mathbb Q^*}(t,\boldsymbol X_t)\Bigr)\dt+Z_t\dW_t\\
 &=-\Bigl(\partial_{x_{k+1}}H(t,\boldsymbol X_t,Y_t)+\tfrac12\partial_{x_{k+1}}a(t,\boldsymbol X_t):\partial_{x_{k+1}}^2v_k^{\mathbb Q^*}(t,\boldsymbol X_t)\\
 &\qquad{}+\partial_{x_{k+1}}f_t(\mathbb Q_{\Pi(t)}^*,\boldsymbol X_t)\Bigr)\dt+Z_t\dW_t,
 \end{aligned}
\]
where the second equality follows from differentiating \eqref{eq:hjb}, and the $i$-th entry of $a:\partial_{x_{k+1}}^3v_k^{\mathbb Q^*}$ is $a:\partial_{x_{k+1}}^2\partial_{x_{k+1}^i}v_k^{\mathbb Q^*}$, with $x_{k+1}^i$ the $i$-th coordinate of $x_{k+1}$, and, for $z\in\R^{d\times d}$ and $1\le j\le k+1$, $\partial_{x_j}a(t,\boldsymbol x):z$ denotes the gradient in $x_j$ of $\boldsymbol x\mapsto a(t,\boldsymbol x):z$. Since $a=\sigma\sigma^\top$ and $\partial_{x_{k+1}}^2v_k^{\mathbb Q^*}$ is symmetric, we have $\tfrac12\partial_{x_j}a(t,\boldsymbol X_t):\partial_{x_{k+1}}^2v_k^{\mathbb Q^*}(t,\boldsymbol X_t)=\partial_{x_j}\sigma(t,\boldsymbol X_t):Z_t$. Hence, taking $j=k+1$,
\[
 \dd Y_t=-\bigl[\partial_{x_{k+1}}H(t,\boldsymbol X_t,Y_t)+\partial_{x_{k+1}}\sigma(t,\boldsymbol X_t):Z_t+\partial_{x_{k+1}}f_t(\mathbb Q_{\Pi(t)}^*,\boldsymbol X_t)\bigr]\dt+Z_t\dW_t.
\]
Since $\mathbb Q_{\Pi(t)}^*=\Law(\boldsymbol X_t)$, this gives the backward equation in~\eqref{eq:fbsde}.

Fix $1\le j<m$. For $t\in[s_k,s_{k+1})$ and $j\le k<m$, we define $M_t^j:=\partial_{x_j}v_k^{\mathbb Q^*}(t,\boldsymbol X_t)$. By the product rule, differentiating~\eqref{eq:hjb} in $x_j$ gives
\[
 \begin{aligned}
 &\partial_t\partial_{x_j}v_k^{\mathbb Q^*}+\partial_{x_{k+1}}\partial_{x_j}v_k^{\mathbb Q^*}\,\partial_yH\bigl(t,\boldsymbol x,\partial_{x_{k+1}}v_k^{\mathbb Q^*}\bigr)\\
 &\quad{}+\tfrac12a(t,\boldsymbol x):\partial_{x_{k+1}}^2\partial_{x_j}v_k^{\mathbb Q^*}{}+\tfrac12\partial_{x_j}a(t,\boldsymbol x):\partial_{x_{k+1}}^2v_k^{\mathbb Q^*}\\
 &\quad{}+\partial_{x_j}H\bigl(t,\boldsymbol x,\partial_{x_{k+1}}v_k^{\mathbb Q^*}\bigr)+\partial_{x_j}f_t(\mathbb Q_{\Pi(t)}^*,\boldsymbol x)=0,
 \end{aligned}
\]
where $\partial_{x_{k+1}}\partial_{x_j}v_k^{\mathbb Q^*}$ is the matrix with $(l,i)$-entry $\partial_{x_{k+1}^i}\partial_{x_j^l}v_k^{\mathbb Q^*}$ and the $l$-th entry of $a:\partial_{x_{k+1}}^2\partial_{x_j}v_k^{\mathbb Q^*}$ is $a:\partial_{x_{k+1}}^2\partial_{x_j^l}v_k^{\mathbb Q^*}$. Since the drift of $X$ is $\partial_yH(t,\boldsymbol X_t,Y_t)$, It\^o's formula gives 
\[
 \begin{aligned}\dd M_t^j& = \; -\bigl[\partial_{x_j}H(t,\boldsymbol X_t,Y_t)+\partial_{x_j}f_t(\mathbb Q_{\Pi(t)}^*,\boldsymbol X_t)
  +\partial_{x_j}\sigma(t,\boldsymbol X_t):Z_t\bigr]\dt\\
 & \qquad +\partial_{x_{k+1}}\partial_{x_j}v_k^{\mathbb Q^*}(t,\boldsymbol X_t)\sigma(t,\boldsymbol X_t)\dW_t,\qquad t\in(s_k,s_{k+1}).\end{aligned}
\]
By \eqref{eq:transmission} and \eqref{eq:terminal}, $M^j$ is continuous on $[s_j,T]$ with $M_T^j=\partial_{x_j}g(\mathbb Q_\Pi^*,\boldsymbol X_T)$. Moreover, by Lemma~\ref{lem:frozenproblem} and Assumption~\ref{ass:coefficients}, both the drift and the diffusion coefficient of $M^j$ are bounded, so the stochastic integral is a martingale. Integrating from $t$ to $T$ and taking conditional expectations given $\F_t$ yields
\[
 \begin{aligned}
 M_t^j={}&\E\Bigl[\partial_{x_j}g(\mathbb Q_\Pi^*,\boldsymbol X_T)+\int_t^T\bigl(\partial_{x_j}H(r,\boldsymbol X_r,Y_r)+\partial_{x_j}f_r(\mathbb Q_{\Pi(r)}^*,\boldsymbol X_r)\\
 &\qquad{}+\partial_{x_j}\sigma(r,\boldsymbol X_r):Z_r\bigr)\dd r\Bigm|\F_t\Bigr],\qquad s_j\le t\le T.
 \end{aligned}
\]
Differentiating \eqref{eq:transmission} in $x$ gives
\[
 \partial_xv_{j-1}^{\mathbb Q^*}(s_j-,\boldsymbol h_{j-1},x)
 =\partial_{h_j}v_j^{\mathbb Q^*}\bigl(s_j,(\boldsymbol h_{j-1},x),x\bigr)+\partial_xv_j^{\mathbb Q^*}\bigl(s_j,(\boldsymbol h_{j-1},x),x\bigr),
\]
Thus, $Y_{s_j-}=M^j_{s_j}+Y_{s_j}$. Since $\mathbb Q_\Pi^*=\Law(\boldsymbol X_T)$ and $\mathbb Q_{\Pi(r)}^*=\Law(\boldsymbol X_r)$, the jump condition in \eqref{eq:fbsde} holds. The terminal condition follows from \eqref{eq:terminal} and $\mathbb Q_\Pi^*=\Law(\boldsymbol X_T)$. Lemma~\ref{lem:frozenproblem} and the boundedness of $\sigma$ give the bounds on $Y$ and $Z$, and Lemma~\ref{lem:state} gives $\E\sup_{t\le T}|X_t|^2<\infty$. Finally, by~\eqref{eq:feedbackbounds}, $y^*$ is bounded by $\bar C$ and Lipschitz continuous in $\boldsymbol x$, uniformly in $t$, so Assumption~\ref{ass:coefficients}\ref{ass:uhat} and~\eqref{eq:uhatbounds} imply that $\hat u(t,\boldsymbol x,y^*(t,\boldsymbol x))$ is Lipschitz and has linear growth in $\boldsymbol x$, uniformly in $t$.
\end{proof}

\subsection{Uniqueness of the FBSDE solution}
\begin{lemma}\label{lem:uniqueness}
Suppose Assumption~\ref{ass:coefficients} holds, and let $(\bar X,\bar Y,\bar Z)$ be either a square-integrable solution of~\eqref{eq:fbsde} with bounded $\bar Y$, or a square-integrable solution of~\eqref{eq:fbsde} when Assumption~\ref{ass:global} holds. Then $(\bar X,\bar Y,\bar Z)=(X,Y,Z)$, where $(X,Y,Z)$ is defined in Lemma~\ref{lem:construction}.
\end{lemma}
\begin{proof}
We define $u_t:=\hat u(t,\bar{\boldsymbol X}_t,\bar Y_t)$ and $\mathbb Q:=\Law(\bar X)$ on $C([0,T];\R^d)$.  Since $\E\sup_{t\le T}|\bar X_t|^2<\infty$, we have $\mathbb Q\in\cP_2(C([0,T];\R^d))$, so Lemma~\ref{lem:frozenproblem} holds for $\mathbb Q$, with $\mathbb Q_{\Pi(t)}=\Law(\bar{\boldsymbol X}_t)$ and $\mathbb Q_\Pi=\Law(\bar{\boldsymbol X}_T)$.

For $1\le j<m$, we denote by $N^j$ a continuous version of 
\[
 \begin{aligned}
 &\E\Bigl[\partial_{x_j}g(\mathbb Q_\Pi,\bar{\boldsymbol X}_T){}+\int_t^T\bigl(\partial_{x_j}H(r,\bar{\boldsymbol X}_r,\bar Y_r)+\partial_{x_j}f_r(\mathbb Q_{\Pi(r)},\bar{\boldsymbol X}_r) \\
 & \qquad \qquad + \; \partial_{x_j}\sigma(r,\bar{\boldsymbol X}_r):\bar Z_r\bigr)\dd r\Bigm|\F_t\Bigr],\quad s_j\le t\le T,
 \end{aligned}
\]
where the random variable inside the conditional expectation is integrable by Assumption~\ref{ass:coefficients}\ref{ass:sigma}, \ref{ass:penalties} and~\ref{ass:bernstein}, since $\E\sup_{r\le T}|\bar Y_r|^2<\infty$ and $\bar Z\in\mathbb H^2(\R^{d\times d})$, so that $\bar Y_{s_j}-\bar Y_{s_j-}=-N_{s_j}^j$ by~\eqref{eq:fbsde}. By the martingale representation theorem, there is a progressively measurable process $\Lambda^j$ with $\int_0^T\|\Lambda_t^j\|^2\dt<\infty$ almost surely such that \[
 \begin{aligned}
 \dd N_t^j={}&-\bigl(\partial_{x_j}H(t,\bar{\boldsymbol X}_t,\bar Y_t)+\partial_{x_j}f_t(\mathbb Q_{\Pi(t)},\bar{\boldsymbol X}_t)\\
 &\quad{}+\partial_{x_j}\sigma(t,\bar{\boldsymbol X}_t):\bar Z_t\bigr)\dt+\Lambda_t^j\dW_t.
 \end{aligned}
\]

For $0\le k<m$ and $t\in[s_k,s_{k+1})$, we compare $\bar Y$, $\bar Z$ and $N^j$, $1\le j\le k$, with
\[
 \tilde Y_t:=\partial_{x_{k+1}}v_k^{\mathbb Q}(t,\bar{\boldsymbol X}_t),\qquad
 \tilde Z_t:=\partial_{x_{k+1}}^2v_k^{\mathbb Q}(t,\bar{\boldsymbol X}_t)\sigma(t,\bar{\boldsymbol X}_t),\qquad
 \tilde M_t^j:=\partial_{x_j}v_k^{\mathbb Q}(t,\bar{\boldsymbol X}_t).
\]
The computation of $\dd Y$ and $\dd M^j$ in the proof of Lemma~\ref{lem:construction} applies to $\tilde Y$ and $\tilde M^j$, except that the drift of $\bar X$ is $\partial_yH(t,\bar{\boldsymbol X}_t,\bar Y_t)$ by~\eqref{eq:envelope}. Hence
\[
 \begin{aligned}
 \dd\tilde Y_t={}&-\bigl[\partial_{x_{k+1}}H(t,\bar{\boldsymbol X}_t,\tilde Y_t)+\partial_{x_{k+1}}\sigma(t,\bar{\boldsymbol X}_t):\tilde Z_t\\
 &\quad{}+\partial_{x_{k+1}}f_t(\mathbb Q_{\Pi(t)},\bar{\boldsymbol X}_t)\bigr]\dt+\tilde Z_t\dW_t\\
 &+\partial_{x_{k+1}}^2v_k^{\mathbb Q}(t,\bar{\boldsymbol X}_t)\bigl(\partial_yH(t,\bar{\boldsymbol X}_t,\bar Y_t)-\partial_yH(t,\bar{\boldsymbol X}_t,\tilde Y_t)\bigr)\dt,\\
 \dd\tilde M_t^j={}&-\bigl[\partial_{x_j}H(t,\bar{\boldsymbol X}_t,\tilde Y_t)+\partial_{x_j}f_t(\mathbb Q_{\Pi(t)},\bar{\boldsymbol X}_t)\\
 &\quad{}+\partial_{x_j}\sigma(t,\bar{\boldsymbol X}_t):\tilde Z_t\bigr]\dt\\
 &+{}\partial_{x_{k+1}}\partial_{x_j}v_k^{\mathbb Q}(t,\bar{\boldsymbol X}_t)\sigma(t,\bar{\boldsymbol X}_t)\dW_t\\
 &+\partial_{x_{k+1}}\partial_{x_j}v_k^{\mathbb Q}(t,\bar{\boldsymbol X}_t)\bigl(\partial_yH(t,\bar{\boldsymbol X}_t,\bar Y_t)-\partial_yH(t,\bar{\boldsymbol X}_t,\tilde Y_t)\bigr)\dt.
 \end{aligned}
\]
Subtracting these equations from the backward equation for $\bar Y$ in~\eqref{eq:fbsde}. Note that $\Law(\bar{\boldsymbol X}_t)=\mathbb Q_{\Pi(t)}$. With the notation $\delta Y:=\bar Y-\tilde Y$, $\delta Z:=\bar Z-\tilde Z$ and $\delta N^j:=N^j-\tilde M^j$, $1\le j\le k$, we have, on $[s_k,s_{k+1})$,
\begin{equation}\label{eq:homogeneous}
 \begin{aligned}
 \dd\delta Y_t=&-\Bigl[\partial_{x_{k+1}}^2v_k^{\mathbb Q}(t,\bar{\boldsymbol X}_t)\,\mathcal C_t\delta Y_t+\mathcal V_t^{k+1}\delta Y_t +\partial_{x_{k+1}}\sigma(t,\bar{\boldsymbol X}_t):\delta Z_t\Bigr]\dt+\delta Z_t\dW_t,\\
 \dd\delta N_t^j=&-\bigl(\partial_{x_{k+1}}\partial_{x_j}v_k^{\mathbb Q}(t,\bar{\boldsymbol X}_t)\,\mathcal C_t+\mathcal V_t^j\bigr)\delta Y_t\dt -\partial_{x_j}\sigma(t,\bar{\boldsymbol X}_t):\delta Z_t\dt\\
 &+\bigl(\Lambda_t^j-\partial_{x_{k+1}}\partial_{x_j}v_k^{\mathbb Q}(t,\bar{\boldsymbol X}_t)\sigma(t,\bar{\boldsymbol X}_t)\bigr)\dW_t,
 \end{aligned}
\end{equation}
where, for $1\le j\le k+1$,
\[
 \mathcal C_t:=\int_0^1\partial_{yy}H\bigl(t,\bar{\boldsymbol X}_t,(1-\lambda)\tilde Y_t+\lambda\bar Y_t\bigr)\dd\lambda,\qquad
 \mathcal V_t^{j}:=\int_0^1\partial_{x_jy}H\bigl(t,\bar{\boldsymbol X}_t,(1-\lambda)\tilde Y_t+\lambda\bar Y_t\bigr)\dd\lambda,
\]
and we used the equality $\mathcal C_t\delta Y_t=\partial_yH(t,\bar{\boldsymbol X}_t,\bar Y_t)-\partial_yH(t,\bar{\boldsymbol X}_t,\tilde Y_t)$ and $\mathcal V_t^{j}\delta Y_t=\partial_{x_j}H(t,\bar{\boldsymbol X}_t,\bar Y_t)-\partial_{x_j}H(t,\bar{\boldsymbol X}_t,\tilde Y_t)$.

By \eqref{eq:feedbackbounds}, $|\tilde Y_t|\le\bar C$, and the derivatives of $v_k^{\mathbb Q}$ above are bounded. The processes $\mathcal C$ and $\mathcal V^{j}$ are also bounded in either case: if $\bar Y$ is bounded, then $(1-\lambda)\tilde Y_t+\lambda\bar Y_t$ is bounded and Assumption~\ref{ass:coefficients}\ref{ass:Hreg} bounds $\mathcal C$ and $\mathcal V^{j}$; if $(\bar X,\bar Y,\bar Z)$ is a square-integrable solution and Assumption~\ref{ass:global} holds, then~\eqref{eq:globalbounds} bounds $\mathcal C$ and $\mathcal V^{j}$. Hence \eqref{eq:homogeneous} is a linear BSDE with bounded coefficients, and $\delta Y$ and $\delta Z$ are square integrable.

By~\eqref{eq:fbsde} and~\eqref{eq:terminal}, $\delta Y_T=0$ and $\delta N_T^j=0$. As in the proof of Lemma~\ref{lem:construction}, differentiating~\eqref{eq:transmission} gives $\tilde Y_{s_j-}=\tilde M_{s_j}^j+\tilde Y_{s_j}$ and shows that $\tilde M^i$ is continuous at $s_j$ for $i<j$. Together with the jump condition for $\bar Y$ and the continuity of $N^i$, this gives
\begin{equation}\label{eq:deltajump}
 \delta Y_{s_j-}=\delta N_{s_j}^j+\delta Y_{s_j},\qquad \delta N_{s_j-}^i=\delta N_{s_j}^i\quad(i<j).
\end{equation}

We prove by backward induction on $k=m-1,\ldots,0$ that $\delta Y=0$, $\delta Z=0$ and $\delta N^i=0$ ($1\le i\le k$) on $[s_k,s_{k+1})$. For $k=m-1$, we have $\delta Y_T=0$ and $\delta N_T^i=0$ ($1\le i<m$). Uniqueness for the square-integrable linear BSDE~\eqref{eq:homogeneous} on $[s_{m-1},T]$ gives $\delta Y=0$ and $\delta Z=0$. Consequently, the drift of $\delta N^i$ is $0$ on $[s_{m-1},T]$. Since $\bar Y=\tilde Y$ and $\bar Z=\tilde Z$ are bounded, $N^i$ and $\delta N^i$ are bounded on $[s_{m-1},T]$, so each $\delta N^i$ ($1\le i<m$) is a bounded martingale on $[s_{m-1},T]$ with $\delta N_T^i=0$, and hence $\delta N^i=0$ on $[s_{m-1},T]$. Now let $1\le k<m$, and suppose that the claim holds for $k$. In particular, $\delta Y_{s_k}=0$ and $\delta N_{s_k}^i=0$ for $1\le i\le k$, so~\eqref{eq:deltajump} with $j=k$ gives $\delta Y_{s_k-}=0$ and $\delta N_{s_k-}^i=0$ for $1\le i<k$. Repeating the same argument on $[s_{k-1},s_k)$ with these terminal values gives $\delta Y=0$, $\delta Z=0$ and $\delta N^i=0$ ($1\le i<k$) on $[s_{k-1},s_k)$. By induction, $\delta Y=0$ and $\delta Z=0$ on $[0,T]$.

Thus $\bar Y_t=\partial_{x_{k+1}}v_k^{\mathbb Q}(t,\bar{\boldsymbol X}_t)$ for $\dd t\otimes\dd\mathbb P$-almost every $(t,\omega)$ with $t\in[s_k,s_{k+1})$, so $u_t=\hat u(t,\bar{\boldsymbol X}_t,\partial_{x_{k+1}}v_k^{\mathbb Q}(t,\bar{\boldsymbol X}_t))$ and $\bar X$ solves \eqref{eq:frozenstate}. Strong uniqueness for \eqref{eq:frozenstate} gives $\bar X=X^{\mathbb Q}$, so $u=u^{\mathbb Q}$ minimizes $J_{\mathbb Q}$ by Lemma~\ref{lem:frozenproblem}. In particular, by~\eqref{eq:uhatbounds} and~\eqref{eq:feedbackbounds}, $|u_t|\le C(1+|\bar{\boldsymbol X}_t|)$ and $u\in\mathbb H^2(\mathbb A)$. Lemma~\ref{lem:state} gives $\bar X=X^u$. Since $\mathbb Q=\Law(X^u)$, the last two terms in~\eqref{eq:linearization} vanish for $u$, and we obtain for every $u'\in\mathbb H^2(\mathbb A)$
\[
 J(u')\ge J_{\mathbb Q}(u')+c_{\mathbb Q}\ge J_{\mathbb Q}(u)+c_{\mathbb Q}=J(u).
\]
Hence $u$ is optimal for~\eqref{eq:cost}, and Lemma~\ref{lem:optimal} gives $u=u^*$. Pathwise uniqueness for~\eqref{eq:state} gives $\bar X=X^{u^*}=X$, so $\mathbb Q=\Law(X)=\Law(X^{\mathbb Q^*})=\mathcal T\mathbb Q^*=\mathbb Q^*$. Since $\mathbb Q=\mathbb Q^*$, $\delta Y=0$ gives $\bar Y=Y$ for $\dd t\otimes\dd\mathbb P$-almost every $(t,\omega)$. Since $\bar Y$ and $Y$ are c\`adl\`ag and have the same terminal value $\partial_{x_m}g(\mathbb Q_\Pi^*,\boldsymbol X_T)$, they are indistinguishable. Finally, $\delta Z=0$ gives $\bar Z=Z$ in $\mathbb H^2(\R^{d\times d})$.
\end{proof}

\subsection{Global well-posedness of the FBSDE}\label{sec:mainproof}
Combining Lemmas~\ref{lem:optimal}, \ref{lem:construction} and~\ref{lem:uniqueness} proves the global well-posedness of the control problem~\eqref{eq:cost} and the McKean--Vlasov FBSDE~\eqref{eq:fbsde}. Specifically, Lemma~\ref{lem:optimal} gives the unique optimal control $u^*$, Lemma~\ref{lem:construction} constructs the square-integrable solution $(X,Y,Z)$ with bounded $Y$ and $Z$ and the feedback representation~\eqref{eq:feedback}, and Lemma~\ref{lem:uniqueness} establishes uniqueness among square-integrable solutions with bounded $Y$ and, under Assumption~\ref{ass:global}, among square-integrable solutions. We have thus proved the following theorem.

\begin{theorem}\label{thm:main}
Let $T>0$ be arbitrary, and suppose Assumption~\ref{ass:coefficients} holds. Then the following hold.
\begin{enumerate}[label=(\alph*)]
\item\label{thm:main-a} The control problem \eqref{eq:cost} has a unique optimal control $u^*\in\mathbb H^2(\mathbb A)$, and the McKean--Vlasov FBSDE \eqref{eq:fbsde} has a square-integrable solution $(X,Y,Z)$ with bounded $Y$ and $Z$ such that $X=X^{u^*}$ and
\begin{equation}\label{eq:feedback}
 u_t^*=\hat u(t,\boldsymbol X_t,Y_t)=\hat u\bigl(t,\boldsymbol X_t,y^*(t,\boldsymbol X_t)\bigr),\qquad t\in[0,T),
\end{equation}
where $y^*(t,\boldsymbol x)\in\R^d$, defined for $t\in[0,T)$ and $\boldsymbol x\in(\R^d)^{\delta^\Pi(t)+1}$, is bounded and Lipschitz continuous in $\boldsymbol x$, uniformly in $t$. The feedback $(t,\boldsymbol x)\mapsto\hat u(t,\boldsymbol x,y^*(t,\boldsymbol x))$ is Lipschitz and has linear growth in $\boldsymbol x$, uniformly in $ t$. 
\item\label{thm:main-b} The solution in \ref{thm:main-a} is the unique square-integrable solution of \eqref{eq:fbsde} with bounded $Y$.
\item\label{thm:main-c} If, in addition, Assumption~\ref{ass:global} holds, then the solution in \ref{thm:main-a} is the unique square-integrable solution of \eqref{eq:fbsde}.
\end{enumerate}
\end{theorem}

Theorem~\ref{thm:main} holds for every horizon $T>0$, every number $m$ of observation times and every value of the constants $L$ and $C_0$ in Assumption~\ref{ass:coefficients}, hence for all nonnegative penalty weights. Section~\ref{sec:examples} applies it to the soft path-law constraints of El Boustany et al.~\cite{elbou2026} and drifts affine in the control, including the soft-constrained Schr\"odinger bridge of Ma et al.~\cite{ma2026softbridge}. Theorem~\ref{thm:main} also identifies the solution of the McKean--Vlasov FBSDE~\eqref{eq:fbsde} with the unique optimal control of~\eqref{eq:cost}, without the joint convexity required by~\cite[Theorem~4.7]{hafsi2026}. By~\eqref{eq:feedback}, the optimal feedback has at most linear growth and depends on the past only through the observations recorded so far.

\section{Applications}\label{sec:examples}
The framework of Section~\ref{sec:model} covers the soft path-law constraints of El Boustany et al.~\cite{elbou2026} (Section~\ref{sec:softpath}) and drifts affine in the control (Section~\ref{sec:affine}), including the soft-constrained Schr\"odinger bridge of Ma et al.~\cite{ma2026softbridge}.

\subsection{Soft path-law constraints and MMD penalties}\label{sec:softpath}

We consider $\mathbb A=\R^d$, $b(t,\boldsymbol x,u)=u$ and $\ell(t,\boldsymbol x,u)=\frac12|u|^2$, with general penalties $F_t$ and $G$, where $F_t$ depends only on the law of the current state, and a diffusion coefficient $\sigma$ that depends only on the current state. This is the soft path-law constraint problem of El Boustany et al.~\cite{elbou2026}. The squared MMD penalties of~\cite{elbou2026} are a special case, which we now introduce.

Let $k$ be a bounded continuous symmetric positive definite kernel on $\R^d\times\R^d$. We denote by $\mathcal H_k$ its reproducing kernel Hilbert space, that is, the completion of the span of $\{k(x,\cdot):x\in\R^d\}$ with respect to the inner product defined by $\langle k(x,\cdot),k(z,\cdot)\rangle_{\mathcal H_k}:=k(x,z)$. Its elements are functions on $\R^d$, and $v(x)=\langle v,k(x,\cdot)\rangle_{\mathcal H_k}$ for all $v\in\mathcal H_k$ and $x\in\R^d$. Since $x\mapsto k(x,\cdot)$ is bounded and continuous from $\R^d$ to $\mathcal H_k$, the reproducing property gives, for probability measures $\mu,\rho$ on $\R^d$,
\[
 \left\langle v,\int k(x,\cdot)(\mu-\rho)(\dd x)\right\rangle_{\mathcal H_k}
 =\int v\,\dd(\mu-\rho),\qquad v\in\mathcal H_k.
\]
Hence the MMD~\cite{gretton2012} satisfies
\begin{align*}
 \operatorname{MMD}_k(\mu,\rho)
 &:=\sup_{\|v\|_{\mathcal H_k}\le1}\left|\int v\,\dd(\mu-\rho)\right|
 =\left\|\int k(x,\cdot)(\mu-\rho)(\dd x)\right\|_{\mathcal H_k},\\
 \operatorname{MMD}_k^2(\mu,\rho)
 &=\iint k(x,z)(\mu-\rho)(\dd x)(\mu-\rho)(\dd z).
\end{align*}
The same construction applies to a kernel $\kappa$ on $\R^{md}\times\R^{md}$. For instance, we take
\begin{equation}\label{eq:mmdpenalties}
 F_t(\mu):=\frac{\lambda_f}{2}\operatorname{MMD}_k^2(\mu,\rho_t),\qquad
 G(\eta):=\frac{\lambda_g}{2}\operatorname{MMD}_\kappa^2(\eta,\rho_\Pi),
 \qquad\lambda_f,\lambda_g\ge0,
\end{equation}
with target laws $\rho_t\in\cP_2(\R^d)$ and $\rho_\Pi\in\cP_2(\R^{md})$.

The following assumption on the kernels and the target laws is a sufficient condition for Assumption~\ref{ass:coefficients}\ref{ass:penalties}; Lemma~\ref{lem:mmd} verifies Assumption~\ref{ass:coefficients}\ref{ass:penalties} under it.

\begin{assumption}\label{ass:mmd}
 The kernels $k$ and $\kappa$ are symmetric, positive definite and of class $C_b^{4+\alpha}$, with $\alpha$ as in Assumption~\ref{ass:coefficients}, and there exists $C\ge0$ such that
\begin{equation}\label{eq:targetregularity}
 d_{\BL}(\rho_t,\rho_s)\le C|t-s|^{1/2},\qquad s,t\in[0,T].
\end{equation}
\end{assumption}

Assumption~\ref{ass:mmd} holds, for instance, for Gaussian kernels with fixed positive bandwidths, their finite nonnegative combinations, and their compositions with feature maps whose derivatives of orders one to five are bounded. 
Examples include the full-path and terminal-state projections in the three-route experiment of~\citet{elbou2026}. The time regularity condition holds for target laws $\rho_t=\Law(\zeta_t)$ of a process $\zeta$ with $\E|\zeta_t-\zeta_s|^2\le C|t-s|$.

\begin{lemma}\label{lem:mmd}
Suppose Assumption~\ref{ass:mmd} holds. Then, after increasing $L$ if necessary, the penalties \eqref{eq:mmdpenalties} satisfy Assumption~\ref{ass:coefficients}\ref{ass:penalties} with the flat derivatives
\[
 f_t(\mu,x)=\lambda_f\int k(x,z)(\mu-\rho_t)(\dd z),\qquad
 g(\eta,\boldsymbol x)=\lambda_g\int\kappa(\boldsymbol x,\boldsymbol z)(\eta-\rho_\Pi)(\dd\boldsymbol z),
\]
and the Lions derivatives of $F_t$ and $G$ are given by $\partial_xf_t(\mu,\cdot)$ and $\partial_{\boldsymbol x}g(\eta,\cdot)$, where
\[
 \begin{aligned}
 \partial_xf_t(\mu,x)&=\lambda_f\int\partial_x k(x,z)(\mu-\rho_t)(\dd z),\\
 \partial_{x_j}g(\eta,\boldsymbol x)&=\lambda_g\int\partial_{x_j}\kappa(\boldsymbol x,\boldsymbol z)(\eta-\rho_\Pi)(\dd\boldsymbol z),\quad 1\le j\le m.
 \end{aligned}
\]
Moreover, $\partial_xf_t$ and $\partial_{\boldsymbol x}g$ are bounded, Lipschitz continuous in the state and, with respect to $d_{\BL}$, in the measure argument, and $\partial_xf_t$ is $1/2$-H\"older continuous in time.
\end{lemma}
\begin{proof}
Fix $\mu,\nu,\rho\in\cP_2(\R^d)$ and $r\in[0,1]$. Since $\operatorname{MMD}_k^2(\mu,\rho)=\iint k(x,z)(\mu-\rho)(\dd x)(\mu-\rho)(\dd z)$, $k$ is symmetric and $(1-r)\mu+r\nu-\rho=(\mu-\rho)+r(\nu-\mu)$, we have
\[
 \begin{aligned}
 \frac12\operatorname{MMD}_k^2\bigl((1-r)\mu+r\nu,\rho\bigr)
 &{}=\frac12\operatorname{MMD}_k^2(\mu,\rho)+r\int\!\left[\int k(x,z)(\mu-\rho)(\dd z)\right](\nu-\mu)(\dd x)\\
 &\quad+\frac{r^2}2\operatorname{MMD}_k^2(\nu,\mu).
 \end{aligned}
\]
Differentiating in $r$ gives
\[
 \frac{\dd}{\dd r}\,\frac12\operatorname{MMD}_k^2\bigl((1-r)\mu+r\nu,\rho\bigr)
 =\int\!\left[\int k(x,z)\bigl((1-r)\mu+r\nu-\rho\bigr)(\dd z)\right](\nu-\mu)(\dd x),
\]
and integrating over $r\in[0,1]$ gives \eqref{eq:flat} for $\frac12\operatorname{MMD}_k^2(\cdot,\rho)$, with the bounded flat derivative $(\mu,x)\mapsto\int k(x,z)(\mu-\rho)(\dd z)$. Moreover, note that
\[
 (1-r)\operatorname{MMD}_k^2(\mu,\rho)+r\operatorname{MMD}_k^2(\nu,\rho)
 =\operatorname{MMD}_k^2\bigl((1-r)\mu+r\nu,\rho\bigr)+r(1-r)\operatorname{MMD}_k^2(\mu,\nu),
\]
so $\operatorname{MMD}_k^2(\cdot,\rho)$ is convex in the sense of Definition~\ref{def:convex}. The same arguments apply to $\kappa$.

 Differentiating the flat derivatives in $x$ and $\boldsymbol x$ gives the formulas for $\partial_xf_t$ and $\partial_{x_j}g$. Since the first derivatives of $k$ and $\kappa$ are bounded and continuous, $\partial_xf_t(\mu,x)$ and $\partial_{\boldsymbol x}g(\eta,\boldsymbol x)$ are bounded and jointly continuous in the state and the measure. Hence $\partial_xf_t(\mu,\cdot)$ and $\partial_{\boldsymbol x}g(\eta,\cdot)$ are the Lions derivatives of $F_t$ at $\mu$ and of $G$ at $\eta$ by~\cite[Proposition~5.48]{carmonadelarue2018}.

The bounds and~\eqref{eq:flatcont} in Assumption~\ref{ass:coefficients}\ref{ass:penalties}, the continuity of $t\mapsto F_t(\delta_0)$ and the properties of $\partial_xf_t$ and $\partial_{\boldsymbol x}g$ follow directly from Assumption~\ref{ass:mmd} and the definition of $d_{\BL}$. This proves Assumption~\ref{ass:coefficients}\ref{ass:penalties}.
\end{proof}

\begin{remark}\label{rem:kl}
The relative entropy $\operatorname{KL}(\mu\,\|\,\rho):=\int\log\frac{\dd\mu}{\dd\rho}\,\dd\mu$, $\mu\ll\rho$, is also convex along mixtures with respect to $ \mu$ in the sense of Definition \ref{def:convex}. However, the flat derivative of $\operatorname{KL}(\mu\,\|\,\rho_t)$ is $\log\frac{\dd\mu}{\dd\rho_t}$, up to an additive constant, which is neither bounded nor continuous in $\mu$, so Assumption~\ref{ass:coefficients}\ref{ass:penalties} fails. A modified Kullback--Leibler (KL) penalty, computed from smooth bins of $\mu$ and $\rho_t$, satisfies Assumption~\ref{ass:coefficients}\ref{ass:penalties}; see Example~\ref{ex:bins} in Appendix~\ref{fv:sec:examples}.
\end{remark}

\begin{corollary}\label{cor:softpath}
Let $\mathbb A=\R^d$, $b(t,\boldsymbol x,u)=u$ and $\ell(t,\boldsymbol x,u)=\frac12|u|^2$.
 Suppose that $\sigma$ depends only on the current state and satisfies Assumption~\ref{ass:coefficients}\ref{ass:sigma}, and $F_t$ and $G$ satisfy Assumption~\ref{ass:coefficients}\ref{ass:penalties}, where $F_t$ depends only on the law of the current state. Then Assumptions~\ref{ass:coefficients} and~\ref{ass:global} hold with $\hat u(t,\boldsymbol x,y)=-y$, $H(t,\boldsymbol x,y)=-\frac12|y|^2$. Consequently, the control problem~\eqref{eq:cost} has a unique optimal control $u^*\in\mathbb H^2(\R^d)$, and the McKean--Vlasov FBSDE~\eqref{eq:fbsde}, which reads
\[
 \begin{cases}
 \dd X_t=-Y_t\dt+\sigma(t,X_t)\dW_t,&X_0=\xi,\\
 \dd Y_t=-\bigl[\partial_x\sigma(t,X_t):Z_t\\
 \qquad\qquad{}+\partial_xf_t(\Law(X_t),X_t)\bigr]\dt+Z_t\dW_t,&t\notin\{s_1,\ldots,s_{m-1}\},\\
 Y_{s_j}-Y_{s_j-}=-\E[\partial_{x_j}g(\Law(\boldsymbol X_T),\boldsymbol X_T)\mid\F_{s_j}],&1\le j<m,\\
 Y_T=\partial_{x_m}g(\Law(\boldsymbol X_T),\boldsymbol X_T),
 \end{cases}
\]
has a unique square-integrable solution $(X,Y,Z)$. They are related by $X=X^{u^*}$ and
\[
 u_t^*=-Y_t=-y^*(t,\boldsymbol X_t),\qquad t\in[0,T),
\]
where $y^*(t,\boldsymbol x)\in\R^d$, defined for $t\in[0,T)$ and $\boldsymbol x\in(\R^d)^{\delta^\Pi(t)+1}$, is bounded and Lipschitz continuous in $\boldsymbol x$, uniformly in $t$. Moreover, $Y$ and $Z$ are bounded $\dd t\otimes\dd\mathbb P$-almost everywhere.
\end{corollary}
\begin{proof}
Conditions~\ref{ass:sigma} and~\ref{ass:penalties} of Assumption~\ref{ass:coefficients} hold by hypothesis, and~\ref{ass:b0} holds since $b(t,\boldsymbol x,0)=0$. For~\ref{ass:drift}, we have $|b(t,\boldsymbol x,u)-b(t,\boldsymbol x,0)|=|u|$, $|b(t,\boldsymbol x,u)-b(t,\boldsymbol x',u)|=0$ and $0\le\ell(t,\boldsymbol x,u)=\frac12|u|^2$. The infimum in~\eqref{eq:reducedhamiltonian} is uniquely attained at $\hat u(t,\boldsymbol x,y)=-y$, which gives~\ref{ass:uhat}, and $H(t,\boldsymbol x,y)=-\frac12|y|^2$, so $H(t,\boldsymbol x,0)=0$ and the derivatives of $H$ in $(\boldsymbol x,y)$ of order 2 to 4 are constant, which gives~\ref{ass:Hreg}. Moreover, $\partial_{\boldsymbol x}H(t,\boldsymbol x,y)=0$, so~\ref{ass:bernstein} holds. For Assumption~\ref{ass:global}, $\partial_{yy}H=-I_d$, $\partial_{\boldsymbol xy}H=0$ and $|\hat u(t,\boldsymbol x,y)|=|y|$, so~\ref{ass:drift} and Assumption~\ref{ass:global} hold for $C_0=\sqrt d$. Substituting $\partial_{\boldsymbol x}H(t,\boldsymbol x,y)=0$, $\partial_{x_j}f_t=0$ and $\partial_{x_j}\sigma=0$, $j\le k$, into~\eqref{eq:fbsde} and applying Theorem~\ref{thm:main}\ref{thm:main-a}, \ref{thm:main-c} gives the conclusion.
\end{proof}

Corollary~\ref{cor:softpath} holds for every horizon $T>0$, every number $m$ of observation times and every value of the constant $L$ in Assumption~\ref{ass:coefficients}\ref{ass:sigma} and~\ref{ass:penalties}, hence for all penalty weights. By Lemma~\ref{lem:mmd}, it applies in particular to the squared MMD with kernels satisfying Assumption~\ref{ass:mmd}. For this penalty, Hafsi et al.~\cite[Theorem~5.15]{hafsi2026} and El Boustany et al.~\cite[Theorem~1]{elbou2026} obtain solvability by a contraction argument, for $T$ below a threshold depending on the weights and the number of observation times. Corollary~\ref{cor:softpath} also identifies the solution of the FBSDE~\eqref{eq:fbsde} with the unique optimal control of~\eqref{eq:cost}, without the convexity of $x\mapsto\sigma(t,x):Z_t$, and its bound on $Y$ shows that the optimal control $u^*=-Y$ is bounded. By the feedback form in Corollary~\ref{cor:softpath}, the optimal control depends on the past only through the observations recorded so far.

For this model, $\hat u(t,\boldsymbol x,y)=-y$, so by Theorem~\ref{thm:main}, the optimal control is characterized by the condition $u^*=-Y$. In practice, a numerical method only produces an approximation $u$ of $u^*$, and the errors $\|u-u^*\|_2$ and $X^u-X$ cannot be computed, since $u^*$ is unknown. We now show that the cost grows quadratically around $u^*$, and that the cost gap $J(u)-J(u^*)$ is bounded by the residual of the optimality condition at $u$. This residual, called the Pontryagin residual, can be computed from $u$ alone and bounds the errors of $u$ (Proposition~\ref{fv:prop:gradient}). It involves the adjoint process of the control $u$, which we now define.
For $u\in\mathbb H^2(\R^d)$, we denote by $(Y^u,Z^u)$ the solution of the linear BSDE with jumps
\begin{equation}\label{fv:adjoint}
 \begin{cases}
 \dd Y_t^u=-\bigl[\partial_xf_t(\Law(X_t^u),X_t^u)\\
 \qquad\quad{}+\partial_x\sigma(t,X_t^u):Z_t^u\bigr]\dt+Z_t^u\dW_t,
       &t\notin\{s_1,\ldots,s_{m-1}\},\\
 Y_{s_j}^u-Y_{s_j-}^u
 =-\E[\partial_{x_j}g(\Law(\boldsymbol X_T^u),\boldsymbol X_T^u)\mid\F_{s_j}],
       &1\le j<m,\\
 Y_T^u=\partial_{x_m}g(\Law(\boldsymbol X_T^u),\boldsymbol X_T^u).
 \end{cases}
\end{equation}
Since its coefficients, jumps and terminal value are bounded, \eqref{fv:adjoint}, solved backwards over the observation intervals, has a unique solution such that $Y^u$ is c\`adl\`ag with $\E\sup_{t\le T}|Y_t^u|^2<\infty$ and $Z^u\in\mathbb H^2(\R^{d\times d})$. By Theorem~\ref{thm:main}, $(Y^{u^*},Z^{u^*})=(Y,Z)$.

\begin{proposition}\label{fv:prop:gradient}
Let $\mathbb A=\R^d$, $b(t,\boldsymbol x,u)=u$ and $\ell(t,\boldsymbol x,u)=\frac12|u|^2$, and suppose that $\sigma$ depends only on the current state and satisfies Assumption~\ref{ass:coefficients}\ref{ass:sigma}, and $F_t$ and $G$ satisfy Assumption~\ref{ass:coefficients}\ref{ass:penalties}, where $F_t$ depends only on the law of the current state. Then for all $u,v\in\mathbb H^2(\R^d)$, the directional derivative $DJ(u)[v]:=\frac{\dd}{\dd\eps}J(u+\eps v)|_{\eps=0}$ exists and
\[
 DJ(u)[v]=\E\int_0^T\langle u_t+Y_t^u,v_t\rangle\dt.
\]
Moreover, there exists $C\ge0$, depending only on $d$, $m$, $T$, the observation times and the constants in Assumption~\ref{ass:coefficients}, such that for all $u\in\mathbb H^2(\R^d)$,
\begin{equation}\label{fv:gradient}
 0\le J(u)-J(u^*)\le C\|u+Y^u\|_2^2,\qquad
 \E\sup_{t\le T}|X_t^u-X_t|^2+\|u-u^*\|_2^2
 \le C\bigl(J(u)-J(u^*)\bigr).
\end{equation}
\end{proposition}
\begin{proof}
Fix $u\in\mathbb H^2(\R^d)$ and first take a bounded direction $v\in\mathbb H^2(\R^d)$.
The derivative of $X^u$ in the direction $v$ is the solution $V$ of
\[
 \dd V_t=v_t\dt+\sum_{r=1}^d
       \partial_x\sigma_{\cdot r}(t,X_t^u)V_t\dd W_t^r,
 \qquad V_0=0.
\]
We apply the product formula to $\langle Y_t^u,V_t\rangle$ on each observation interval. The cross variation of the diffusion terms cancels the term $\partial_x\sigma(t,X_t^u):Z_t^u$ in \eqref{fv:adjoint}, and the jump at $s_j$ is $\langle Y_{s_j}^u-Y_{s_j-}^u,V_{s_j}\rangle$, whose expectation is $-\E\langle\partial_{x_j}g(\Law(\boldsymbol X_T^u),\boldsymbol X_T^u),V_{s_j}\rangle$ since $V_{s_j}$ is $\F_{s_j}$-measurable. Summing over the intervals gives
\[
 DJ(u)[v]=\E\int_0^T\langle u_t+Y_t^u,v_t\rangle\dt.
\]
The bounded derivatives in Assumption~\ref{ass:coefficients} justify the differentiation for bounded $v$, and the formula extends to $v\in\mathbb H^2$ by approximation.

We next prove the second inequality in \eqref{fv:gradient}. Let $\mathbb Q^*=\Law(X)$ be the fixed point of the map $\mathcal T$ in~\eqref{eq:bestresponse}. For this model, $b(t,\boldsymbol x,u)\cdot y+\ell(t,\boldsymbol x,u)-H(t,\boldsymbol x,y)=\frac12|u+y|^2$, so the integrand in~\eqref{eq:verification} is $\frac12|u_t+\partial_{x_{k+1}}v_k^{\mathbb Q}(t,\boldsymbol X_t^u)|^2$. By \eqref{eq:linearization} with $\mathbb Q=\mathbb Q^*$ and \eqref{eq:verification},
\[
 J(u)-J(u^*)\ge J_{\mathbb Q^*}(u)-J_{\mathbb Q^*}(u^*)=\frac12
 \sum_{k=0}^{m-1}\E\int_{s_k}^{s_{k+1}}|u_t+\partial_{x_{k+1}}v_k^{\mathbb Q^*}(t,\boldsymbol X_t^u)|^2\dt.
\]
Since $u_t^*=-\partial_{x_{k+1}}v_k^{\mathbb Q^*}(t,\boldsymbol X_t^{u^*})$ for $t\in[s_k,s_{k+1})$, the feedback is Lipschitz continuous in $(\boldsymbol h_k,x)$ by \eqref{eq:feedbackbounds}, and $|\boldsymbol X_t^u-\boldsymbol X_t^{u^*}|^2\le(m+1)\sup_{s\le t}|X_s^u-X_s^{u^*}|^2$, the Burkholder--Davis--Gundy and Gronwall inequalities, applied successively on the observation intervals, yield the second inequality in \eqref{fv:gradient}.

For the first inequality, let $\mathbb Q:=\Law(X^u)$ and $r:=u+Y^u$. By \eqref{eq:linearization}, whose last two terms vanish for $u$ and are nonnegative for $u^*$, and by \eqref{eq:verification},
\[
 J(u)-J(u^*)\le J_{\mathbb Q}(u)-J_{\mathbb Q}(u^*)
 \le\frac12\sum_{k=0}^{m-1}\E\int_{s_k}^{s_{k+1}}|u_t+\partial_{x_{k+1}}v_k^{\mathbb Q}(t,\boldsymbol X_t^u)|^2\dt.
\]
For $1\le j<m$, we denote by $N^{u,j}$ a continuous version of the martingale $\E[\partial_{x_j}g(\Law(\boldsymbol X_T^u),\boldsymbol X_T^u)\mid\F_t]$, so that $Y_{s_j}^u-Y_{s_j-}^u=-N_{s_j}^{u,j}$, and by $\Lambda^j$ the integrand in its martingale representation.    
For $0\le k<m$ and $t\in[s_k,s_{k+1})$, we compare $Y^u$, $Z^u$ and $N^{u,j}$, $1\le j\le k$, with
\[
 \tilde Y_t:=\partial_{x_{k+1}}v_k^{\mathbb Q}(t,\boldsymbol X_t^u),\qquad
 \tilde Z_t:=\partial_{x_{k+1}}^2v_k^{\mathbb Q}(t,\boldsymbol X_t^u)\sigma(t,X_t^u),\qquad
 \tilde M_t^j:=\partial_{x_j}v_k^{\mathbb Q}(t,\boldsymbol X_t^u),
\]
and write $\delta Y:=Y^u-\tilde Y$, $\delta Z:=Z^u-\tilde Z$ and $\delta N^j:=N^{u,j}-\tilde M^j$. Then $u_t+\partial_{x_{k+1}}v_k^{\mathbb Q}(t,\boldsymbol X_t^u)=r_t-\delta Y_t$. The computation in the proof of Lemma~\ref{lem:uniqueness}, with the drift $u_t$ of $X^u$ in place of $\partial_yH(t,\bar{\boldsymbol X}_t,\bar Y_t)$, gives on $[s_k,s_{k+1})$
\[
 \begin{aligned}
 \dd\delta Y_t&=-\bigl[\partial_{x_{k+1}}^2v_k^{\mathbb Q}(t,\boldsymbol X_t^u)(r_t-\delta Y_t)+\partial_x\sigma(t,X_t^u):\delta Z_t\bigr]\dt+\delta Z_t\dW_t,\\
 \dd\delta N_t^j&=-\partial_{x_{k+1}}\partial_{x_j}v_k^{\mathbb Q}(t,\boldsymbol X_t^u)(r_t-\delta Y_t)\dt+\bigl(\Lambda_t^j-\partial_{x_{k+1}}\partial_{x_j}v_k^{\mathbb Q}(t,\boldsymbol X_t^u)\sigma(t,X_t^u)\bigr)\dW_t.
 \end{aligned}
\]
Moreover, $\delta Y_T=0$, $\delta N_T^j=0$, and \eqref{eq:deltajump} holds as in the proof of Lemma~\ref{lem:uniqueness}. By~\eqref{eq:feedbackbounds} and Assumption~\ref{ass:coefficients}\ref{ass:sigma}, the coefficients are bounded. Applying It\^o's formula to $|\delta Y_t|^2$ and $|\delta N_t^j|^2$ on each $[s_k,s_{k+1})$ together with Young's and Gronwall's inequalities gives
\[
 \sup_{t\in[s_k,s_{k+1})}\Bigl(\E|\delta Y_t|^2+\sum_{j=1}^k\E|\delta N_t^j|^2\Bigr)
 \le C\Bigl(\E|\delta Y_{s_{k+1}-}|^2+\sum_{j=1}^k\E|\delta N_{s_{k+1}-}^j|^2+\E\int_{s_k}^{s_{k+1}}|r_s|^2\ds\Bigr).
\]
Iterating backward from $k=m-1$ to $0$ using $\delta Y_T=0$, $\delta N_T^j=0$ and~\eqref{eq:deltajump} yields $\|\delta Y\|_2\le C\|r\|_2$, and hence $J(u)-J(u^*)\le\frac12\|r-\delta Y\|_2^2\le C\|r\|_2^2$, proving the first inequality in~\eqref{fv:gradient}.
\end{proof}

\subsection{Drifts affine in the control}\label{sec:affine}
\begin{corollary}\label{cor:affine}
Suppose that $\sigma$ satisfies Assumption~\ref{ass:coefficients}\ref{ass:sigma}, $F_t$ and $G$ satisfy Assumption~\ref{ass:coefficients}\ref{ass:penalties}, and $\mathbb A=\R^q$. Assume that $b$ is affine in $u$, $b(t,\boldsymbol x,u)=\bar b(t,\boldsymbol x)+B(t,\boldsymbol x)u$, where $\bar b$ satisfies the conditions on $b(\cdot,\cdot,0)$ in Assumption~\ref{ass:coefficients}\ref{ass:b0}, and $B(t,\boldsymbol x)\in\R^{d\times q}$ and its spatial derivatives up to order 4 are bounded and continuous. Assume further that $\ell$ is continuous and such that, for every $R>0$, its derivatives in $(\boldsymbol x,u)$ up to order 4 are continuous and bounded on $\{|u|\le R\}$. Finally, assume that there exist $C_0\ge0$ and $c_0>0$ such that
\[
 \begin{gathered}
 \partial_{uu}\ell(t,\boldsymbol x,u)\ge c_0I_q,\qquad
 -C_0\le\ell(t,\boldsymbol x,u)\le C_0(1+|u|^2),\\
 |\partial_{\boldsymbol x}\ell(t,\boldsymbol x,u)|\le C_0\bigl(1+|\ell(t,\boldsymbol x,u)|\bigr).
 \end{gathered}
\]
\begin{enumerate}[label=(\alph*)]
\item\label{cor:affine-a} Then Assumption~\ref{ass:coefficients} holds with a larger $C_0$, and Theorem~\ref{thm:main}\ref{thm:main-a}, \ref{thm:main-b} apply.
\item\label{cor:affine-b} If, in addition, $B$ does not depend on $\boldsymbol x$ and $\partial_{\boldsymbol x}\partial_u\ell$ is bounded, then Assumption~\ref{ass:global} holds, after enlarging $C_0$, and Theorem~\ref{thm:main}\ref{thm:main-c} applies.
\end{enumerate}
\end{corollary}
\begin{proof}
Throughout this proof, $C$ denotes a generic constant independent of $t$, $\boldsymbol x$, $y$ and $u$, which may change from line to line.

(a) Since $\partial_{uu}\ell\ge c_0I_q$, for every $(t,\boldsymbol x,y)$ the function $u\mapsto b(t,\boldsymbol x,u)\cdot y+\ell(t,\boldsymbol x,u)$ is strictly convex and tends to infinity as $|u|\to\infty$. Hence it has a unique minimizer $\hat u(t,\boldsymbol x,y)$, which is the unique zero in $\R^q$ of $u\mapsto B(t,\boldsymbol x)^\top y+\partial_u\ell(t,\boldsymbol x,u)$. 
By $\partial_{uu}\ell\ge c_0I_q$ and this first-order condition $B(t,\boldsymbol x)^\top y+\partial_u\ell(t,\boldsymbol x,\hat u)=0$, $\hat u=\hat u(t,\boldsymbol x,y)$ satisfies
\[
 c_0|\hat u|^2
 \le\bigl(\partial_u\ell(t,\boldsymbol x,\hat u)-\partial_u\ell(t,\boldsymbol x,0)\bigr)\cdot \hat u
 =-\bigl(B(t,\boldsymbol x)^\top y+\partial_u\ell(t,\boldsymbol x,0)\bigr)\cdot \hat u.
\]
Since $B$ and $\partial_u\ell(\cdot,\cdot,0)$ are bounded, this yields $|\hat u(t,\boldsymbol x,y)|\le C(1+|y|)$. The implicit function theorem shows that $\hat u$ is three times continuously differentiable in $(\boldsymbol x,y)$, with derivatives that are bounded and continuous on $\{|y|\le R\}$ for every $R>0$. By the chain rule and the first-order condition $B(t,\boldsymbol x)^\top y+\partial_u\ell(t,\boldsymbol x,\hat u)=0$, we have
\[
 \begin{aligned}
 H(t,\boldsymbol x,y)&=\bigl(\bar b(t,\boldsymbol x)+B(t,\boldsymbol x)\hat u\bigr)\cdot y+\ell(t,\boldsymbol x,\hat u),\\
 \partial_yH(t,\boldsymbol x,y)&=\bar b(t,\boldsymbol x)+B(t,\boldsymbol x)\hat u,\\
 \partial_{\boldsymbol x}H(t,\boldsymbol x,y)&=\partial_{\boldsymbol x}b(t,\boldsymbol x,\hat u)^\top y+\partial_{\boldsymbol x}\ell(t,\boldsymbol x,\hat u),
 \end{aligned}
\]
and $\partial_{\boldsymbol x}b(t,\boldsymbol x,u)-\partial_{\boldsymbol x}\bar b(t,\boldsymbol x)$ is linear in $u$, with coefficients given by the spatial derivatives of $B$. Together with the regularity of $\hat u$, $\ell$, $B$ and $\bar b$, these formulas show that the derivatives of $\partial_{\boldsymbol x}H$ and $\partial_yH$ in $(\boldsymbol x,y)$ of order 1 to 3 are bounded and continuous on $\{|y|\le R\}$. Since $H(t,\boldsymbol x,0)=\ell(t,\boldsymbol x,\hat u(t,\boldsymbol x,0))\in[-C_0,C_0]$, this gives~\ref{ass:Hreg}; the regularity of $\hat u$ gives~\ref{ass:uhat}, and \ref{ass:sigma}--\ref{ass:drift} follow from the assumptions.

For~\ref{ass:bernstein}, the formula for $\partial_{\boldsymbol x}H$, the bounds on $\ell$ and $\partial_{\boldsymbol x}\ell$, and $|\hat u(t,\boldsymbol x,y)|\le C(1+|y|)$ give $|\partial_{\boldsymbol x}H(t,\boldsymbol x,y)|\le C(1+|\hat u|)|y|+C_0\bigl(1+C_0(1+|\hat u|^2)\bigr)\le C(1+|y|^2)$.

(b) Here $\hat u(t,\boldsymbol x,y)$ is the unique zero of $B(t)^\top y+\partial_u\ell(t,\boldsymbol x,u)$ in $\R^q$, and the implicit function theorem gives
\[
 \partial_y\hat u=-(\partial_{uu}\ell)^{-1}B(t)^\top,\qquad
 \partial_{\boldsymbol x}\hat u=-(\partial_{uu}\ell)^{-1}\partial_{\boldsymbol x}\partial_u\ell,
\]
which are bounded by $c_0^{-1}\|B\|_\infty$ and $c_0^{-1}\|\partial_{\boldsymbol x}\partial_u\ell\|_\infty$. Since $\partial_yH=\bar b(t,\boldsymbol x)+B(t)\hat u$ and $\partial_{\boldsymbol x}\bar b$ is bounded, the matrices $\partial_{yy}H$ and $\partial_{\boldsymbol xy}H$ are bounded.
\end{proof}

Under our regularity assumptions, taking $m=1$, $F_t=0$, $\ell(t,\boldsymbol x,u)=\frac12|u|^2$ and $B=\sigma(t)$ in Corollary~\ref{cor:affine} recovers the soft-constrained Schr\"odinger bridge model of~\citet{ma2026softbridge}. The corollary also applies when $B=\sigma(t,x)$. It gives uniqueness in the full square-integrable class by~\ref{cor:affine-b} when $\sigma$ depends only on time, and in the square-integrable class with bounded $Y$ by~\ref{cor:affine-a} when $\sigma$ depends on the state. 
 Corollary~\ref{cor:affine}\ref{cor:affine-a} also covers running costs quadratic in $u$ with a state-dependent, uniformly positive definite matrix that, together with its spatial derivatives up to order 4, is bounded and continuous.
 
Assumption~\ref{ass:coefficients} also holds for the running cost $(1+|u|^2)^{3/4}$ and $B=I_d$, since then $\partial_{\boldsymbol x}H(t,\boldsymbol x,y)=\partial_{\boldsymbol x}\bar b(t,\boldsymbol x)^\top y$. For this cost, $H(t,\boldsymbol x,y)-\bar b(t,\boldsymbol x)\cdot y$ behaves like $-c|y|^3$ as $|y|\to\infty$ for some $c>0$. This lies beyond the theory of BSDEs with quadratic growth~\cite{kobylanski2000}, but Lemma~\ref{lem:hjb} still holds.

\section{Conclusion}\label{sec:conclusion}


Theorem~\ref{thm:main} establishes the existence of a unique optimal
control in Lipschitz feedback form for every finite horizon and all
nonnegative penalty weights. The drift, diffusion coefficient and
running cost may depend on the current state and on previously recorded
states. The associated McKean--Vlasov FBSDE with jumps is well posed in
the class of square-integrable solutions with bounded adjoint $Y$.
Under Assumption~\ref{ass:global}, uniqueness holds in the full
square-integrable class.

The proof uses classical solutions of nonlinear HJB equations on
successive observation intervals, with the recorded states as
parameters. Differentiating the matching conditions yields the jumps
of the adjoint. Uniform estimates for the HJB equations allow us to
apply the Schauder--Tychonoff theorem to the best-response map on path
laws. At the fixed point, convexity of the penalties along mixtures of
laws and the HJB verification identity imply optimality for the
McKean--Vlasov control problem.

When the control is the drift, the running control cost is quadratic,
the diffusion depends on the current state, and the running law
penalties act on its marginals, Proposition~\ref{fv:prop:gradient}
combines quadratic growth and gradient domination. 
The cost gap, the squared control error and the mean-square supremum error of the
state are bounded by $C\|u+Y^u\|_{\mathbb H^2}^2$, with $C$ independent
of $u$. 
This gives an a posteriori error bound for approximate
controls in the soft path-law formulation of time-series generation.

\appendix
\section{Proofs of the lemmas}\label{app:proofs}

\subsection{Proof of Lemma~\ref{lem:state}}\label{app:state}
\begin{proof}
On each interval $[s_k,s_{k+1})$, the past observations are
$\mathcal F_{s_k}$-measurable. The drift is Lipschitz in the current
state with coefficient $C_0(1+|u_t|)$, which is integrable in time
almost surely, and $\sigma$ is bounded and uniformly Lipschitz.
Since $|b(t,\boldsymbol x,u_t)|\le C(1+|\boldsymbol x|+|u_t|)$,
the monotonicity and growth conditions of \cite{gyongyKrylov1980}
hold on each interval. Applying their result successively gives a
unique strong solution on $[0,T]$. The Burkholder--Davis--Gundy and Gronwall inequalities yield
$\E\sup_{t\le T}|X_t^u|^2<\infty$.
\end{proof}

\subsection{Proof of Lemma~\ref{lem:hjb}}\label{app:hjb}
Throughout this subsection, $C$ denotes a generic constant that may change from line to line, and $C_R$ a constant that may also depend on $R$. We prove Lemma~\ref{lem:hjb} by 3 steps: we first establish uniqueness and global $L^\infty$ and gradient a priori bounds via a  maximum principle (Lemma~\ref{lem:hjb-apriori}), next construct the global $C_b^4$-in-$x$ classical solution (Lemmas~\ref{lem:hjb-splitting} and~\ref{lem:hjb-spatial}), and finally obtain the parameter derivatives in $h$ up to order $4$ by the probabilistic representation (Lemma~\ref{lem:hjb-param}).

\begin{lemma}\label{lem:hjb-apriori}
The following hold for~\eqref{eq:paramhjb}.
\begin{enumerate}[label=(\alph*)]
\item\label{lem:hjb-apriori-a} Equation~\eqref{eq:paramhjb} has at most one bounded classical solution on $[\underline t,\bar t]\times\R^p\times\R^d$ with $\partial_xv$ bounded.
\item\label{lem:hjb-apriori-b} There exist constants $K,M>0$, independent of the starting time $s$,
the parameter $h$, and the particular solution $v$, such that for every $s\in[\underline t,\bar t)$ and every function $v$ that solves~\eqref{eq:paramhjb} classically with $v\in
C_b^{0,0,3}([s,\bar t]\times\R^p\times\R^d), v(\bar t,\cdot,\cdot)=\Phi$, we have
\begin{equation}\label{eq:gradbound}
 \|v\|_\infty<K,\qquad \|\partial_xv\|_\infty\le M.
\end{equation}
\end{enumerate}
\end{lemma}
\begin{proof}
Fix $s\in[\underline t,\bar t)$ and $h\in\R^p$, and omit $h$ from the notation throughout this proof. Consider a linear operator 
$$\mathcal A (\bar b)f:=\partial_tf+\frac12a:\partial_{xx}f+\bar b\cdot\partial_xf $$ 
with drift satisfying $|\bar b(t,x)|\le C(1+|x|)$ for some constant $C$ (which may depend on $h$). Let $\mathcal B(t,x):=e^{C'(\bar t-t)}(1+|x|^2)$ for $(t,x) \in [s,\bar t)\times\R^d$. Since $a$ is bounded, it satisfies $\mathcal A(\bar b)\mathcal B\le-\mathcal B<0$ once $C'$ is chosen large enough. 

Let $ g$ be a bounded function. Let $w$ be bounded and continuous on $[s,\bar t]\times\R^d$, $C^{1,2}$ on $[s,\bar t)\times\R^d$, and solves $\mathcal A(\bar b)w+g=0$. For each $n\ge1$, consider the function 
$$ w_n:=w-(\bar t-t)\|g\|_\infty-\mathcal B/n. $$
Since $ w_n$ is continuous and tends to \(-\infty\) as \(|x|\to\infty\),
uniformly in \(t\in[s,\bar t]\), it attains its maximum at some
\((t_n,x_n)\in[s,\bar t]\times\R^d\). If $t_n<\bar t$, then $\partial_xw_n=0$, $\partial_{xx}w_n\le0$ and $\partial_tw_n\le0$ at $(t_n,x_n)$. This contradicts $\mathcal A(\bar b)w_n=-g+\|g\|_\infty-\mathcal A(\bar b)\mathcal B/n\ge\mathcal B/n>0$. Hence $t_n=\bar t$, so $w_n(t,x)\le w(\bar t,x_n)-\mathcal B(\bar t,x_n)/n\le\|w(\bar t,\cdot)\|_\infty$. Letting $n\to\infty$ and repeating this analysis for $-w$, we obtain
\begin{equation}\label{eq:linmax}
 \|w\|_\infty\le\|w(\bar t,\cdot)\|_\infty+(\bar t-s)\|g\|_\infty.
\end{equation}

With the bound above, we prove the desired results in this lemma. 

For~\ref{lem:hjb-apriori-a}, let $v$ and $\tilde v$ be bounded classical solutions of~\eqref{eq:paramhjb} on $[\underline t,\bar t]\times\R^p\times\R^d$ with bounded $\partial_xv$ and $\partial_x\tilde v$. Let $w:=\tilde v-v$. Then $ w$ solves $\mathcal A(\bar b)w=0$ with $w(\bar t,\cdot)=0$ and $\bar b(t,x):=\int_0^1\partial_yH(t,(h,x),\lambda\partial_x\tilde v+(1-\lambda)\partial_xv)\,\dd\lambda$. By \eqref{eq:uhatbounds}, \eqref{eq:envelope} and Assumption~\ref{ass:coefficients}\ref{ass:b0}--\ref{ass:drift}, $|\partial_yH(t,(h,x),y)|\le C_R(1+|h|+|x|)$ on $\{|y|\le R\}$, so $\bar b$ grows at most linearly in $x$ and~\eqref{eq:linmax} yields $\tilde v=v$.

For~\ref{lem:hjb-apriori-b}, let $ v$ be a solution of \eqref{eq:paramhjb} on $ [s,\bar t)\times\R^p\times\R^d$. Then we have 
\[
\mathcal A(\bar b)v+U+H(\cdot,\cdot,0)=0,
\]
where $ \bar b(t,x) =\int_0^1\partial_yH(t,x,\lambda\partial_xv) \dd\lambda$. Since $ \bar b$ has linear growth with respect to $ x$, \eqref{eq:linmax} yields
\[
 \|v\|_\infty<K:=1+\|\Phi\|_\infty+T\bigl(\|U\|_\infty+\|H(\cdot,\cdot,0)\|_\infty\bigr).
\]
We now bound $\partial_xv$ by a Bernstein-type argument (cf.~\cite{barles2021}). 
Define 
$$\mathcal Gf:=\partial_tf+\partial_yH(t,x,\partial_xv)\cdot\partial_xf+\frac12a:\partial_{xx}f.$$
Differentiating \eqref{eq:paramhjb} with respect to $ x_i$, we have
\[
 \mathcal Gv=\partial_yH(t,x,\partial_xv)\cdot\partial_xv-H(t,x,\partial_xv)-U,\qquad
 \mathcal G\partial_{x_i}v=-\tfrac12\partial_{x_i}a:\partial_{xx}v-\partial_{x_i}H(t,x,\partial_xv)-\partial_{x_i}U.
\]
Note that $$\mathcal G\frac12|\partial_xv|^2=\sum_i\partial_{x_i}v\,\mathcal G\partial_{x_i}v+\frac12\operatorname{tr}(\partial_{xx}v\,a\,\partial_{xx}v)\ge-C(1+|\partial_xv|^2)-\partial_xv\cdot\partial_xH(t,x,\partial_xv),$$ 
where we used the uniform ellipticity $a\ge\underline aI_d$, and Young's inequality on $\partial_{x_i}a:\partial_{xx}v$.

By Assumption~\ref{ass:coefficients}\ref{ass:bernstein}, the term $-\partial_xv\cdot\partial_xH(t,x,\partial_xv)$ is only bounded below by $-C(|\partial_xv|+|\partial_xv|^3)$, this lower bound does not ensure positivity for large \( |\partial_xv| \) and the maximum principle cannot be applied directly to $\frac12|\partial_xv|^2$. To create a positive term of order $|\partial_xv|^4$ that dominates these negative terms, we divide $\frac12|\partial_xv|^2$ by a smooth positive function $\mathcal N(v)$ on $[-K,K]$. Setting $\widetilde{\mathcal G}f:=\mathcal Gf+\frac{\mathcal N'(v)}{\mathcal N(v)}\partial_xv^\top a\,\partial_xf$, the product rule gives
\[
 \widetilde{\mathcal G}\frac{\frac12|\partial_xv|^2}{\mathcal N(v)}
 =\frac{\mathcal G(\frac12|\partial_xv|^2)}{\mathcal N(v)}-\frac{\mathcal N'(v)\,\frac12|\partial_xv|^2}{\mathcal N(v)^2}\,\mathcal Gv-\frac{\mathcal N''(v)\,\frac12|\partial_xv|^2}{2\mathcal N(v)^2}\,\partial_xv^\top a\,\partial_xv.
\]
Since $a\ge\underline aI_d$, choosing $\mathcal N$ strictly concave ($\mathcal N''<0$) makes the last term positive of order $|\partial_xv|^4$. By~\eqref{eq:concavity}, we have $
\mathcal Gv\le -H(t,x,0)-U(t,x)
\le \|H(\cdot,\cdot,0)\|_\infty+\|U\|_\infty$. Choosing $\mathcal N$ increasing ($\mathcal N'>0$) bounds the second term from below by $-C|\partial_xv|^2$. We therefore choose
\[
 \mathcal N(r):=1+\frac{r+K}{2K}-\frac{(r+K)^2}{16K^2}\qquad(r\in[-K,K]),
\]
which satisfies $\mathcal N\in[1,\frac74]$, $\mathcal N'\in[\frac1{4K},\frac1{2K}]$ and $\mathcal N''=-\frac1{8K^2}$. With this choice, the last term is at least $\frac{\underline a}{98K^2}|\partial_xv|^4$, the second is at least $-C|\partial_xv|^2$ by~\eqref{eq:concavity}, and Assumption~\ref{ass:coefficients}\ref{ass:bernstein} together with Young's inequality gives $|\partial_xv\cdot\partial_xH(t,x,\partial_xv)|\le\frac{\underline a}{196K^2}|\partial_xv|^4+C(1+|\partial_xv|^2)$. Since $\mathcal N\ge1$, we obtain
\[
 \widetilde{\mathcal G}\frac{\frac12|\partial_xv|^2}{\mathcal N(v)}\ge\frac{\underline a}{196K^2}|\partial_xv|^4-C(1+|\partial_xv|^2),
\]
which is strictly positive wherever $|\partial_xv|>C$, for a constant $C$ independent of $s$ and $h$.

We now apply the maximum principle. Note that
$\widetilde{\mathcal G}=\mathcal A(\bar b_v)$, where
$
\bar b_v(t,x)
:=\partial_yH(t,x,\partial_xv)
  +\frac{\mathcal N'(v)}{\mathcal N(v)}\,a\,\partial_xv.$
Since $\partial_xv$ is bounded, the first term has linear growth with respect to $x$ and
the second term is bounded. Thus $|\bar b_v(t,x)|\le C(1+|x|)$.

The function $ \mathcal{B}$ in the proof of part \ref{lem:hjb-apriori-a} satisfies
$\widetilde{\mathcal G}\mathcal B\le-\mathcal B$
for $C'$ sufficiently large. For each $n\ge1$, the function $\frac12|\partial_xv|^2/\mathcal N(v)-\mathcal B/n$ attains its maximum on $[s,\bar t]\times\R^d$ at some $(t_n,x_n)$. If $t_n<\bar t$, then $\widetilde{\mathcal G}(\frac12|\partial_xv|^2/\mathcal N(v)-\mathcal B/n)(t_n,x_n)\le0$, hence $\widetilde{\mathcal G}(\frac12|\partial_xv|^2/\mathcal N(v))(t_n,x_n)\le\widetilde{\mathcal G}\mathcal B(t_n,x_n)/n<0$, which makes $|\partial_xv(t_n,x_n)|\le C$ for some constant $ C > 0$. On the other hand, if $t_n=\bar t$, then $|\partial_xv(t_n,x_n)|\le\|\partial_x\Phi\|_\infty$. In either case, since $\mathcal N\ge1$,
\[
 \frac{\frac12|\partial_xv(t,x)|^2}{\mathcal N(v(t,x))}-\frac{\mathcal B(t,x)}{n}
 \le\frac{\frac12|\partial_xv(t_n,x_n)|^2}{\mathcal N(v(t_n,x_n))}
 \le\tfrac12\max\{\|\partial_x\Phi\|_\infty,C\}^2.
\]
Letting $n\to\infty$ and using $\mathcal N\le\frac74$ yields $\|\partial_xv\|_\infty\le M:=\frac{\sqrt7}2\max\{\|\partial_x\Phi\|_\infty,C\}$. This proves the desired result.
\end{proof}

To construct a solution with four spatial derivatives by solving a mild form of~\eqref{eq:paramhjb} via contraction, we face two issues: $H(t,\boldsymbol x,y)$ may grow linearly in $\boldsymbol x$ (so we cannot treat all of $H$ as a bounded perturbation of $\frac12a:\partial_{xx}$), and closing a fixed-point map in $C_b^4$ when the nonlinearity depends on $\partial_xv$ requires the linear semigroup to gain one spatial derivative. Rather than absorbing $\partial_yH(\cdot,\cdot,0)$ into the linear drift (which would require fifth derivatives of $H$ when differentiating four times in $\boldsymbol x$), we extract a symmetric difference in $y$.

For the standard basis $e_1,\ldots,e_d$ of $\R^d$, define
\[
 \widetilde b_i(t,\boldsymbol x)
 :=\tfrac12\bigl(H(t,\boldsymbol x,e_i)-H(t,\boldsymbol x,-e_i)\bigr),\qquad
 H_0(t,h,x,y):=H(t,(h,x),y)-\widetilde b(t,(h,x))\cdot y.
\]

\begin{lemma}\label{lem:hjb-splitting}
\begin{enumerate}[label=(\alph*)]
\item\label{lem:hjb-splitting-a}  $\sup_t\bigl(|\widetilde b(t,0)|
+\|\partial_{\boldsymbol x}\widetilde b(t,\cdot)\|_{C_b^3}\bigr)<\infty$, and $\partial_{\boldsymbol x}\widetilde b\in C_b^{0,3}([\underline t,\bar t)\times\R^{p+d})$. In particular,
$|\widetilde b(t,\boldsymbol x)|\le C(1+|\boldsymbol x|)$.

\item\label{lem:hjb-splitting-b} For every $R>0$, let $S_R:=\R^{p+d}\times\{y\in\R^d:|y|\le R\}$.
Then $H_0\in C_b^{0,4}([\underline t,\bar t)\times S_R)$.

\item\label{lem:hjb-splitting-c} For $(t,h,x)\in[\underline t,\bar t]\times\R^p\times\R^d$, denote by $\Theta^{t,h,x}$ the solution on $[t,\bar t]$ of
\[
 \dd\Theta_r=\widetilde b\bigl(r,(h,\Theta_r)\bigr)\dd r+\sigma\bigl(r,(h,\Theta_r)\bigr)\dW_r,\qquad \Theta_t=x.
\]
Then, for every $q\ge1$,
\[
 \sup_{\underline t\le t\le s\le\bar t}\ \sup_{(h,x)\in\R^{p+d}}
 \E\bigl[|\partial_x^\gamma\Theta_s^{t,h,x}|^q\bigr]<\infty
 \qquad(1\le|\gamma|\le4),
\]
and $(t,s,h,x)\mapsto\partial_x^\gamma\Theta_s^{t,h,x}$ is continuous in $L^q$ for $0\le|\gamma|\le4$.
Moreover, for every $\psi\in C_b^j(\R^d),j\in\{0,\ldots,4\}$ and $\underline t\le t<s\le\bar t$, we have
\begin{equation}\label{eq:smoothing}
\begin{aligned}
 \sup_{h\in \R^p}\bigl\|\E\psi(\Theta_s^{t,h,\cdot})\bigr\|_{C_b^{j+1}}
 &\le C(s-t)^{-1/2}\|\psi\|_{C_b^j}\quad(0\le j\le3),\\
 \sup_{h\in \R^p}\bigl\|\E\psi(\Theta_s^{t,h,\cdot})\bigr\|_{C_b^j}&\le C\|\psi\|_{C_b^j}\quad(0\le j\le4).
\end{aligned}
\end{equation}
\end{enumerate}
\end{lemma}
\begin{proof}
Since $\widetilde b_i(t,\boldsymbol x)=\frac12\int_{-1}^1\partial_{y_i}H(t,\boldsymbol x,\lambda e_i)\,\dd\lambda$, Assumption~\ref{ass:coefficients}\ref{ass:Hreg} gives $|\partial_yH(t,\boldsymbol x,y)-\widetilde b(t,\boldsymbol x)|\le C_R$ for $|y|\le R$. By~\eqref{eq:envelope} and Assumption~\ref{ass:coefficients}\ref{ass:b0}--\ref{ass:uhat}, $\partial_yH(t,0,0)=b(t,0,\hat u(t,0,0))$ is bounded, so $\widetilde b(t,0)$ is bounded. By~\ref{ass:Hreg} and~\ref{ass:bernstein}, the spatial derivatives of $\widetilde b$ of orders $1$ to $4$ are bounded and continuous, proving~\ref{lem:hjb-splitting-a}.

For~\ref{lem:hjb-splitting-b}, on $\{|y|\le R\}$, $\partial_yH_0$ and the derivatives of $H_0$ of orders $2$ to $4$ in $(\boldsymbol x,y)$ are bounded. Since $H_0(t,\cdot,0)=H(t,\cdot,0)$ and $\partial_{\boldsymbol x\boldsymbol x}H_0(t,\cdot,0)$ are bounded, Taylor's formula bounds $\partial_{\boldsymbol x}H_0(t,\cdot,0)$; integrating $\partial_yH_0$ and $\partial_{\boldsymbol xy}H_0$ in $y$ then bounds $H_0$ and $\partial_{\boldsymbol x}H_0$ on $\{|y|\le R\}$.

For~\ref{lem:hjb-splitting-c}, the $L^q$-moment bounds and $L^q$-continuity of the $x$-derivatives of $\Theta_s^{t,h,x}$ up to order $4$ follow from Assumption~\ref{ass:coefficients}\ref{ass:sigma} and~\ref{lem:hjb-splitting-a}. By the Bismut--Elworthy--Li formula~\cite{elworthy1994}, for bounded measurable $\psi:\R^d\to\R$, $e\in\R^d$ and $\underline t\le t<s\le\bar t$,
\begin{equation}\label{eq:bel}
 \partial_x\E\bigl[\psi(\Theta_s^{t,h,x})\bigr]e
 =\frac1{s-t}\E\left[\psi(\Theta_s^{t,h,x})
 \int_t^s\left\langle\sigma\bigl(r,(h,\Theta_r^{t,h,x})\bigr)^{-1}
               \partial_x\Theta_r^{t,h,x}e,\dd W_r\right\rangle\right].
\end{equation}
Differentiating $\E[\psi(\Theta_s^{t,h,x})]$ and~\eqref{eq:bel} in $x$ and applying the Burkholder--Davis--Gundy inequality gives~\eqref{eq:smoothing}.
\end{proof}

Combining the result of Lemma~\ref{lem:hjb-splitting} with the global gradient bound of Lemma~\ref{lem:hjb-apriori} yields the unique global classical solution with four spatial derivatives.

\begin{lemma}\label{lem:hjb-spatial}
Equation~\eqref{eq:paramhjb} admits a unique bounded classical solution $v$ on $[\underline t,\bar t]\times\R^p\times\R^d$ with $\partial_xv$ bounded. Moreover, $v\in C_b^{0,0,4}([\underline t,\bar t]\times\R^p\times\R^d)$ and
\begin{equation}\label{eq:apriori}
 \sup_{t\in[\underline t,\bar t],\,h\in\R^p}\|v(t,h,\cdot)\|_{C_b^4}< \infty.
\end{equation}
\end{lemma}
\begin{proof}
For $0<\delta\le\bar t-\underline t$, the space $\mathcal X_\delta :=C_b^{0,0,4}([\bar t-\delta,\bar t]\times\R^p\times\R^d)$ is a Banach space under the norm $\|w\|_{\mathcal X_\delta}:=\sup_{t\in[\bar t-\delta,\bar t],\,h\in\R^p}\|w(t,h,\cdot)\|_{C_b^4}$. Define the mild operator
\[
 (\mathcal Mw)(t,h,x):=\E\Bigl[\Phi\bigl(h,\Theta_{\bar t}^{t,h,x}\bigr)
 +\int_t^{\bar t}\bigl(U(s,h,\Theta_s^{t,h,x})+H_0\bigl(s,h,\Theta_s^{t,h,x},\partial_xw(s,h,\Theta_s^{t,h,x})\bigr)\bigr)\ds\Bigr].
\]
For $w,\widetilde w$ in the closed ball of radius $R\ge1$ of $\mathcal X_\delta$, Lemma~\ref{lem:hjb-splitting}\ref{lem:hjb-splitting-b} and the chain rule give $\|H_0(s,h,\cdot,\partial_xw(s,h,\cdot))\|_{C_b^3}\le C_R$ and
\[
 \|H_0(s,h,\cdot,\partial_xw(s,h,\cdot))-H_0(s,h,\cdot,\partial_x\widetilde w(s,h,\cdot))\|_{C_b^3}\le C_R\|w-\widetilde w\|_{\mathcal X_\delta}.
\]
By~\eqref{eq:smoothing} and the $L^q$-continuity of $\Theta_s^{t,h,x}$ and its $x$-derivatives (handling the integral near $s=t$ via~\eqref{eq:smoothing}), $\mathcal M$ maps this ball into $\mathcal X_\delta$ and satisfies
\begin{equation}\label{eq:mild-contraction}
 \|\mathcal Mw\|_{\mathcal X_\delta}\le C+C_R\sqrt\delta,\qquad
 \|\mathcal Mw-\mathcal M\widetilde w\|_{\mathcal X_\delta}\le C_R\sqrt\delta\,\|w-\widetilde w\|_{\mathcal X_\delta}.
\end{equation}
Choosing $R:=C+1$ and $\delta>0$ with $C_R\sqrt\delta\le\frac12$, the contraction mapping theorem yields a local fixed point $v\in\mathcal X_\delta$ of
\begin{equation}\label{eq:mild}
 v(t,h,x)=\E\Bigl[\Phi\bigl(h,\Theta_{\bar t}^{t,h,x}\bigr)
 +\int_t^{\bar t}\bigl(U(s,h,\Theta_s^{t,h,x})+H_0\bigl(s,h,\Theta_s^{t,h,x},\partial_xv(s,h,\Theta_s^{t,h,x})\bigr)\bigr)\ds\Bigr].
\end{equation}

Any solution $v\in\mathcal X_{\bar t-s}$ of~\eqref{eq:mild} on an interval $[s,\bar t]\subset[\underline t,\bar t]$ solves~\eqref{eq:paramhjb} classically on $[s,\bar t)$. Indeed, for $t\in[s,\bar t)$, fixed $(h,x)$ and $0<\eta<\bar t-t$, the Markov property of $\Theta$ and It\^o's formula for $v(t+\eta,h,\Theta_r)$ along $\Theta_r:=\Theta_r^{t,h,x}$ give
\[
 \begin{aligned}
 v(t,h,x)-&v(t+\eta,h,x)=\E\int_t^{t+\eta}\Bigl(\widetilde b\bigl(r,(h,\Theta_r)\bigr)\cdot\partial_xv(t+\eta,h,\Theta_r)\\
 &+\tfrac12a\bigl(r,(h,\Theta_r)\bigr):\partial_{xx}v(t+\eta,h,\Theta_r)
 +U(r,h,\Theta_r)+H_0\bigl(r,h,\Theta_r,\partial_xv(r,h,\Theta_r)\bigr)\Bigr)\dd r.
 \end{aligned}
\]
Dividing by $\eta$, letting $\eta\downarrow0$ and using $|\widetilde b(r,(h,\Theta_r))|\le C(1+|h|+|\Theta_r|)$, by dominated convergence, the right derivative of $v$ in $t$ exists, is continuous, and equals $-\frac12a:\partial_{xx}v-H(\cdot,\cdot,\partial_xv)-U$, so $v$ is continuously differentiable in $t$ and solves~\eqref{eq:paramhjb} classically on $[s,\bar t)$.

To extend the local solution backwards to $[\underline t,\bar t]$, let $v\in\mathcal X_{\bar t-s}$ solve~\eqref{eq:mild} on $[s,\bar t]$. As shown above, $v$ also solves~\eqref{eq:paramhjb} on $[s,\bar t)$, so Lemma~\ref{lem:hjb-apriori}\ref{lem:hjb-apriori-b} gives $\sup_{t\in[s,\bar t],\,h\in\R^p}\|v(t,h,\cdot)\|_{C_b^1}\le C$ and $|\partial_xv|\le M$. We bound $\|v(t,h,\cdot)\|_{C_b^{j+1}}$ inductively for $j=1,2,3$. Given a uniform $C_b^j$ bound on $v$, Lemma~\ref{lem:hjb-splitting}\ref{lem:hjb-splitting-b} on $\{|y|\le M\}$ and the chain rule yield $\|H_0(r,h,\cdot,\partial_xv(r,h,\cdot))\|_{C_b^j}\le C(1+\|v(r,h,\cdot)\|_{C_b^{j+1}})$. Applying~\eqref{eq:smoothing} to~\eqref{eq:mild} then gives
\[
 \|v(t,h,\cdot)\|_{C_b^{j+1}}
 \le C+C\int_t^{\bar t}(r-t)^{-1/2}\|v(r,h,\cdot)\|_{C_b^{j+1}}\,\dd r
 \qquad(t\in[s,\bar t],\ h\in\R^p).
\]
By a Volterra Gronwall inequality, $\sup_{t\in[s,\bar t],\,h\in\R^p}\|v(t,h,\cdot)\|_{C_b^{j+1}}\le C$ independently of $s$. By induction, we have
\[
 \sup_{t\in[s,\bar t],\,h\in\R^p}\|v(t,h,\cdot)\|_{C_b^4}\le C,
\]
with $C$ independent of $s$. 

Apply~\eqref{eq:mild-contraction} with terminal time $s$ and terminal
condition $v(s,\cdot,\cdot)$. Since $\sup_h\|v(s,h,\cdot)\|_{C_b^4}$
is bounded independently of $s$, we may choose the same radius $R$
and the same $\delta>0$ with $C_R\sqrt{\delta}\le\tfrac12$ at each step.
The contraction mapping theorem then gives a solution on $[s',s]$,
where $s':=\max\{\underline t,s-\delta\}$.

For $t\in[s',s]$, apply~\eqref{eq:mild} on $[t,s]$ with terminal
condition $v(s,h,\cdot)$. Substituting the representation of
$v(s,h,\cdot)$ on $[s,\bar t]$ and using the Markov property of $\Theta$
at time $s$ gives~\eqref{eq:mild} on $[t,\bar t]$ with terminal
condition $\Phi(h,\cdot)$.
Repeating this step finitely many times, we obtain a classical solution
$v\in\mathcal X_{\bar t-\underline t}$ satisfying~\eqref{eq:apriori}.
Uniqueness follows from Lemma~\ref{lem:hjb-apriori}\ref{lem:hjb-apriori-a}.
\end{proof}

\begin{lemma}\label{lem:hjb-param}
Let $v$ be the solution from Lemma~\ref{lem:hjb-spatial}, define $b^*(s,h,x):=\partial_yH(s,(h,x),\partial_xv(s,h,x))$, and for $(t,h,x)\in[\underline t,\bar t]\times\R^p\times\R^d$ denote by $\widehat\Theta=\widehat\Theta^{t,h,x}$ the solution of
\[
 \dd\widehat\Theta_s=b^*(s,h,\widehat\Theta_s)\ds+\sigma\bigl(s,(h,\widehat\Theta_s)\bigr)\dW_s,\qquad \widehat\Theta_t=x,
\]
and set $y_s:=\partial_xv(s,h,\widehat\Theta_s)$. Fix $1\le i\le p$, denote by $\sigma_l$ the $l$-th column of $\sigma$, and let $V=V^{t,h,x}$ solve
\[
 \dd V_s=\sum_{l=1}^d\Bigl(\partial_{h_i}\sigma_l\bigl(s,(h,\widehat\Theta_s)\bigr)+\partial_x\sigma_l\bigl(s,(h,\widehat\Theta_s)\bigr)V_s\Bigr)\dd W_s^l,\qquad V_t=0.
\]
\begin{enumerate}[label=(\alph*)]
\item\label{lem:hjb-param-a} The derivative $\partial_{h_i}v$ exists and satisfies
\begin{equation}\label{eq:firstvariation}
\begin{aligned}
 \partial_{h_i}v(t,h,x)
 &=\E\Bigl[\partial_{h_i}\Phi(h,\widehat\Theta_{\bar t})
             +\partial_x\Phi(h,\widehat\Theta_{\bar t})\cdot V_{\bar t}\\
 &\qquad+\int_t^{\bar t}\bigl(\partial_{h_i}H(s,(h,\widehat\Theta_s),y_s)
               +\partial_xH(s,(h,\widehat\Theta_s),y_s)\cdot V_s\\
 &\qquad\qquad\quad+\partial_{h_i}U(s,h,\widehat\Theta_s)
               +\partial_xU(s,h,\widehat\Theta_s)\cdot V_s\bigr)\ds\Bigr].
\end{aligned}
\end{equation}
\item\label{lem:hjb-param-b} We have $v\in C_b^{0,4}([\underline t,\bar t]\times\R^{p+d})$. Moreover, $\partial_{h_i}v$ is continuously differentiable in $t$ and satisfies the equation obtained by differentiating~\eqref{eq:paramhjb} with respect to $h_i$,
\begin{equation}\label{eq:paramhjb-hi}
 \begin{aligned}
 &\partial_t\partial_{h_i}v+\tfrac12a\bigl(t,(h,x)\bigr):\partial_{xx}\partial_{h_i}v+\partial_yH\bigl(t,(h,x),\partial_xv\bigr)\cdot\partial_x\partial_{h_i}v\\
 &\quad{}+\tfrac12\partial_{h_i}a\bigl(t,(h,x)\bigr):\partial_{xx}v+\partial_{h_i}H\bigl(t,(h,x),\partial_xv\bigr)+\partial_{h_i}U(t,h,x)=0.
 \end{aligned}
\end{equation}
\end{enumerate}
\end{lemma}
\begin{proof}
For~\ref{lem:hjb-param-a}, we first prove that $\partial_xv$ is differentiable with respect to $h$ by comparing $v(t,h',x)$ and $v(t,h,x)$. Consider the process $\widetilde\Theta^{h'}$ whose drift is frozen at $b^*(s,h,\widehat\Theta_s)$ and whose diffusion uses $h'\in\R^p$:
\[
 \dd\widetilde\Theta_s^{h'}=b^*(s,h,\widehat\Theta_s)\ds+\sigma\bigl(s,(h',\widetilde\Theta_s^{h'})\bigr)\dW_s,\qquad \widetilde\Theta_t^{h'}=x.
\]
Since $\sigma$ is Lipschitz in $(h,x)$, the Burkholder--Davis--Gundy and Gronwall inequalities give $\E\sup_{s\in[t,\bar t]}|\widetilde\Theta_s^{h'}-\widehat\Theta_s|^q\le C_q|h'-h|^q$ for every $q\ge2$. Applying It\^o's formula to $v(s,h',\widetilde\Theta_s^{h'})$ on $[t,\bar t-\eta]$ using~\eqref{eq:paramhjb} and letting $\eta\downarrow0$, we have
\[
 v(t,h',x)=\E\Bigl[\Phi(h',\widetilde\Theta_{\bar t}^{h'})+\int_t^{\bar t}\Bigl(U(s,h',\widetilde\Theta_s^{h'})+H\bigl(s,(h',\widetilde\Theta_s^{h'}),y_s'\bigr)-b^*(s,h,\widehat\Theta_s)\cdot y_s'\Bigr)\ds\Bigr],
\]
where $y_s':=\partial_xv(s,h',\widetilde\Theta_s^{h'})$.
Set $h'=h$ and let $ y_s := \partial_xv(s,h,\widetilde\Theta_s^{h})$. We obtain 
\begin{equation}\label{eq:vdiff}
 \begin{aligned}
 v(t,h',x)-v(t,h,x)=\E\Bigl[&\Phi(h',\widetilde\Theta_{\bar t}^{h'})-\Phi(h,\widehat\Theta_{\bar t})\\
 &+\int_t^{\bar t}\Bigl(U(s,h',\widetilde\Theta_s^{h'})-U(s,h,\widehat\Theta_s)\\
 &\qquad+H\bigl(s,(h',\widetilde\Theta_s^{h'}),y_s'\bigr)-H\bigl(s,(h,\widehat\Theta_s),y_s'\bigr)+\mathcal R_s(h')\Bigr)\ds\Bigr],
 \end{aligned}
\end{equation}
where
\[
 \mathcal R_s(h'):=H\bigl(s,(h,\widehat\Theta_s),y_s'\bigr)-H\bigl(s,(h,\widehat\Theta_s),y_s\bigr)-\partial_yH\bigl(s,(h,\widehat\Theta_s),y_s\bigr)\cdot(y_s'-y_s).
\] By concavity of $H(s,(h,\widehat\Theta_s),\cdot)$, we have $\mathcal R_s(h')\le0$. Since $\Phi$, $U$ and $H(\cdot,\cdot,y)$ (on $\{|y|\le M\}$, by Assumption~\ref{ass:coefficients}\ref{ass:bernstein}) are uniformly Lipschitz in $(h,x)$, this yields $v(t,h',x)-v(t,h,x)\le C|h'-h|$. Swapping $h$ and $h'$ gives 
\begin{equation}\label{eq:v-lipschitz-h}
  \sup_{t,x}|v(t,h',x)-v(t,h,x)|\le C|h'-h|.
\end{equation}
Set $h'=h+\eps e_i$ and $w(x):=v(s,h',x)-v(s,h,x)$. By~\eqref{eq:v-lipschitz-h},
$\|w\|_\infty\le C|\eps|$, while~\eqref{eq:apriori} gives
$\|\partial_x^3w\|_\infty\le C$. The interpolation inequality \cite[Lemma A.3.3]{fernandez-real2024IntegroDifferential}
$\|\partial_xw\|_\infty\le
C\|w\|_\infty^{2/3}\|\partial_x^3w\|_\infty^{1/3}$
and the bound on $\partial_{xx}v$ in~\eqref{eq:apriori} yield
\[
 \|\partial_xv(s,h',\cdot)-\partial_xv(s,h,\cdot)\|_\infty\le C|\eps|^{2/3},
 \qquad
 |y_s'-y_s|\le C\bigl(|\eps|^{2/3}+|\widetilde\Theta_s^{h'}-\widehat\Theta_s|\bigr).
\]
By the bound on $\partial_{yy}H$ on $\{|y|\le M\}$, we have $|\mathcal R_s(h')|\le C|y_s'-y_s|^2$, so $|\eps|^{-1}\E\int_t^{\bar t}|\mathcal R_s(h')|\ds\le C|\eps|^{1/3}$. This, along with Taylor's formula on $\sigma$ and Gronwall's inequality, yields
\[
\begin{gathered}
 \lim_{\eps\to0}\frac1\eps\E\int_t^{\bar t}|\mathcal R_s(h+\eps e_i)|\ds=0,\\
 \lim_{\eps\to0}\sup_{s\in[t,\bar t]}\E\left|\frac{\widetilde\Theta_s^{h+\eps e_i}-\widehat\Theta_s}\eps-V_s\right|^2=0.
\end{gathered}
\]
Dividing~\eqref{eq:vdiff} at $h'=h+\eps e_i$ by $\eps$ and letting $\eps\to0$ gives~\eqref{eq:firstvariation}.

For~\ref{lem:hjb-param-b}, we prove by induction on $r=0,1,2,3,4$ that
\[
 \partial_h^\beta\partial_x^\gamma v\in C_b^{0,0}\bigl([\underline t,\bar t]\times\R^{p+d}\bigr)
 \qquad(|\beta|\le r,\ |\beta|+|\gamma|\le4).
\]
The case $r=0$ holds by Lemma~\ref{lem:hjb-spatial}. Fix $r\in\{1,2,3,4\}$ and assume the claim for $r-1$. Because $v$ enters $b^*(s,h,x)=\partial_yH(s,(h,x),\partial_xv(s,h,x))$ only through $\partial_xv$, the chain rule and Assumption~\ref{ass:coefficients}\ref{ass:sigma}, \ref{ass:Hreg} imply that $b^*$ is continuous and that the derivatives $\partial_h^{\tilde\beta}\partial_x^{\tilde\gamma}$ of $b^*$ ($1\le|\tilde\beta|+|\tilde\gamma|\le3$) and of $\sigma,\partial_{h_i}\sigma_l,\partial_x\sigma_l$ ($0\le|\tilde\beta|+|\tilde\gamma|\le3$) are bounded and continuous whenever $|\tilde\beta|\le r-1$. 

Differentiating the equations for $\widehat\Theta$ and $V$ in $(h,x)$
gives the mixed derivatives
$\partial_h^{\beta'}\partial_x^\gamma(\widehat\Theta,V)$.
By dominated convergence and the Burkholder--Davis--Gundy and Gronwall
inequalities, these derivatives are continuous in $(t,s,h,x)$ in
$L^q(\Omega)$ and satisfy, for every $1\le q<\infty$,
\[
 \sup_{\underline t\le t\le s\le\bar t}\ \sup_{(h,x)\in\R^{p+d}}
 \E\bigl[|\partial_h^{\beta'}\partial_x^\gamma\widehat\Theta_s^{t,h,x}|^q+|\partial_h^{\beta'}\partial_x^\gamma V_s^{t,h,x}|^q\bigr]<\infty
\]
for all $|\beta'|\le r-1$ and $1\le|\beta'|+|\gamma|\le3$ (and also for $\beta'=\gamma=0$ in the $V$-term). Thus, taking the derivatives on the both side of \eqref{eq:firstvariation}, we obtain the result for $ r$. By the induction, we prove that $ v\in C_b^{0,4}([\underline t,\bar t]\times\R^{p+d})$.

Finally, writing~\eqref{eq:paramhjb} as $v(t,h,x)=\Phi(h,x)+\int_t^{\bar t}(\frac12a:\partial_{xx}v+H(\cdot,\cdot,\partial_xv)+U)(s,h,x)\ds$ and differentiating under the integral sign in $h_i$ shows that $\partial_{h_i}v$ is continuously differentiable in $t$ and solves \eqref{eq:paramhjb-hi} on $[\underline t,\bar t)\times\R^p\times\R^d$.
\end{proof}

\begin{proof}[Proof of Lemma~\ref{lem:hjb}]
Existence, uniqueness, and boundedness of $v$ and $\partial_xv$ are given by Lemma~\ref{lem:hjb-spatial}, and $v\in C_b^{0,4}([\underline t,\bar t]\times\R^{p+d})$ (together with \eqref{eq:paramhjb-hi}, used in Lemmas~\ref{lem:construction} and~\ref{lem:uniqueness}) follow from Lemma~\ref{lem:hjb-param}\ref{lem:hjb-param-b}.
\end{proof}

\subsection{Proof of Lemma~\ref{lem:frozenproblem}}\label{app:frozen}
\begin{proof}
Throughout this proof, $C$ denotes a generic constant independent of $\mathbb Q$, which may change from line to line.

We first construct the solutions $v_k^{\mathbb Q}$ of~\eqref{eq:hjb}--\eqref{eq:transmission} and establish their uniform $C_b^4$ bounds. For $t\in[s_k,s_{k+1})$, set $U_t:=f_t(\mathbb Q_{\Pi(t)},\cdot)$.
By~\eqref{eq:flatcont} and Assumption~\ref{ass:coefficients}\ref{ass:penalties},
\[
 \|U_t-U_s\|_{C_b^1}
 \le L\bigl(|t-s|^{1/2}+d_{\BL}(\mathbb Q_{\Pi(t)},\mathbb Q_{\Pi(s)})\bigr),
 \qquad
 \|U_t-U_s\|_{C_b^{4+\alpha}}\le 2L.
\]
Hence, with $\theta=\alpha/(3+\alpha)$, the interpolation inequality \cite[Lemma A.3.3]{fernandez-real2024IntegroDifferential} gives
\[
 \|U_t-U_s\|_{C_b^4}
 \le C\|U_t-U_s\|_{C_b^1}^{\theta}
       \|U_t-U_s\|_{C_b^{4+\alpha}}^{1-\theta}
 \le C\bigl(|t-s|^{1/2}
       +d_{\BL}(\mathbb Q_{\Pi(t)},\mathbb Q_{\Pi(s)})\bigr)^{\theta}.
\]
By continuity of the paths, $t\mapsto\mathbb Q_{\Pi(t)}$ is weakly
continuous on $[s_k,s_{k+1})$ and has a weak left limit at $s_{k+1}$.
Thus $t\mapsto U_t$ extends continuously to $[s_k,s_{k+1}]$ in $C_b^4$.
Also, $g(\mathbb Q_\Pi,\cdot)\in C_b^4(\R^{md})$ by
Assumption~\ref{ass:coefficients}\ref{ass:penalties}, and \eqref{eq:transmission} preserves $C_b^4$ regularity across observation times. Applying Lemma~\ref{lem:hjb} backward on $[s_k,s_{k+1}]$ for $k=m-1,\ldots,0$ (with $h=\boldsymbol h_k$) yields the solutions $v_k^{\mathbb Q}$. The derivative bounds are uniform in $\mathbb Q$, since the bounds on
$f_t$ and $g$ in Assumption~\ref{ass:coefficients}\ref{ass:penalties}
are uniform in the measure argument. This proves \eqref{eq:feedbackbounds}.

Next, we establish the well-posedness of~\eqref{eq:frozenstate} and the admissibility of $u^{\mathbb Q}$. By~\eqref{eq:envelope}, the drift in~\eqref{eq:frozenstate} is $\partial_yH(t,\boldsymbol x,\partial_{x_{k+1}}v_k^{\mathbb Q}(t,\boldsymbol x))$ at $\boldsymbol x=\boldsymbol X_t^{\mathbb Q}$.  By~\eqref{eq:feedbackbounds} and Assumption~\ref{ass:coefficients}\ref{ass:b0}--\ref{ass:Hreg}, this drift is Lipschitz in $\boldsymbol x$ and has at most linear growth, uniformly in $t$ and $\mathbb Q$.  Together with Assumption~\ref{ass:coefficients}\ref{ass:sigma}, this ensures that \eqref{eq:frozenstate} has a unique strong solution with $\E\sup_{t\le T}|X_t^{\mathbb Q}|^2<\infty$.  Hence $u^{\mathbb Q}\in\mathbb H^2(\mathbb A)$ by \eqref{eq:uhatbounds} and \eqref{eq:feedbackbounds}.

We now verify that $u^{\mathbb Q}$ is the unique minimizer of $J_{\mathbb Q}$ in $\mathbb H^2(\mathbb A)$. For $u\in\mathbb H^2(\mathbb A)$, It\^o's formula and \eqref{eq:hjb} on each interval $[s_k,s_{k+1})$ give
\[
 \begin{aligned}
 \dd v_k^{\mathbb Q}(t,\boldsymbol X_t^u)
 ={}&\Bigl(b(t,\boldsymbol X_t^u,u_t)\cdot\partial_{x_{k+1}}v_k^{\mathbb Q}(t,\boldsymbol X_t^u)-f_t(\mathbb Q_{\Pi(t)},\boldsymbol X_t^u)-H\bigl(t,\boldsymbol X_t^u,\partial_{x_{k+1}}v_k^{\mathbb Q}(t,\boldsymbol X_t^u)\bigr)\Bigr)\dt\\
 &+\partial_{x_{k+1}}v_k^{\mathbb Q}(t,\boldsymbol X_t^u)^\top\sigma(t,\boldsymbol X_t^u)\dW_t.
 \end{aligned}
\]
Summing over the intervals and taking expectations, we obtain
\begin{equation}\label{eq:verification}
 \begin{aligned}
 J_{\mathbb Q}(u)
 ={}&\E v_0^{\mathbb Q}(0,\xi)
  +\sum_{k=0}^{m-1}\E\int_{s_k}^{s_{k+1}}\Bigl(b(t,\boldsymbol X_t^u,u_t)\cdot\partial_{x_{k+1}}v_k^{\mathbb Q}(t,\boldsymbol X_t^u)\\
 &\quad+\,\ell(t,\boldsymbol X_t^u,u_t)-H\bigl(t,\boldsymbol X_t^u,\partial_{x_{k+1}}v_k^{\mathbb Q}(t,\boldsymbol X_t^u)\bigr)\Bigr)\dt,
 \end{aligned}
\end{equation}
where the stochastic integrals have zero expectation since $\partial_xv_k^{\mathbb Q}$ and $\sigma$ are bounded, and all terms are integrable by Assumption~\ref{ass:coefficients}\ref{ass:penalties}--\ref{ass:Hreg} and Lemma~\ref{lem:state}. By \eqref{eq:reducedhamiltonian}, the integrand is nonnegative, and by Assumption~\ref{ass:coefficients}\ref{ass:uhat}, it vanishes if and only if $u_t=\hat u(t,\boldsymbol X_t^u,\partial_{x_{k+1}}v_k^{\mathbb Q}(t,\boldsymbol X_t^u))$. Hence $u^{\mathbb Q}$ is admissible and attains the lower bound $\E v_0^{\mathbb Q}(0,\xi)$ in \eqref{eq:verification}. By \eqref{eq:verification}, every minimizer $u$ satisfies this identity for $\dd t\otimes\dd\mathbb P$-almost every $(t,\omega)$ with $t\in[s_k,s_{k+1})$, so $X^u$ solves \eqref{eq:frozenstate}, and strong uniqueness gives $X^u=X^{\mathbb Q}$ and $u=u^{\mathbb Q}$.

Finally, we prove \eqref{eq:frozengrad} under~\eqref{eq:frozenhyp}. By~\eqref{eq:hjb}, $w_k:=v_k^{\mathbb Q^n}-v_k^{\mathbb Q}$ solves
\[
 \partial_tw_k+\tfrac12a:\partial_{xx}w_k+\bar b_k\cdot\partial_xw_k+f_t(\mathbb Q_{\Pi(t)}^n,\cdot)-f_t(\mathbb Q_{\Pi(t)},\cdot)=0
 \qquad\text{on }[s_k,s_{k+1}),
\]
where $\bar b_k(t,\boldsymbol h_k,x):=\int_0^1\partial_yH\bigl(t,(\boldsymbol h_k,x),(1-\lambda)\partial_xv_k^{\mathbb Q}+\lambda\partial_xv_k^{\mathbb Q^n}\bigr)\,\dd\lambda$. By \eqref{eq:uhatbounds}, \eqref{eq:envelope}, \eqref{eq:feedbackbounds} and Assumption~\ref{ass:coefficients}\ref{ass:b0}--\ref{ass:drift}, $\bar b_k$ has at most linear growth in $x$. Applying~\eqref{eq:linmax} on each $[s_k,s_{k+1}]$ and iterating backward via~\eqref{eq:terminal}--\eqref{eq:transmission} yields
\[
 \max_k\|w_k\|_\infty
 =\max_k\|v_k^{\mathbb Q^n}-v_k^{\mathbb Q}\|_\infty
 \le T\sup_{t<T,\boldsymbol x}|f_t(\mathbb Q_{\Pi(t)}^n,\boldsymbol x)-f_t(\mathbb Q_{\Pi(t)},\boldsymbol x)|+\sup_{\boldsymbol x}|g(\mathbb Q_\Pi^n,\boldsymbol x)-g(\mathbb Q_\Pi,\boldsymbol x)|,
\]
which tends to $0$ as $n\to\infty$ by \eqref{eq:frozenhyp}. Since $\|\partial_{xx}w_k\|_\infty\le 2\bar C$ uniformly in $n$ by~\eqref{eq:feedbackbounds}, the interpolation inequality~\cite[Lemma~A.3.3]{fernandez-real2024IntegroDifferential} gives
\[
 \|\partial_xv_k^{\mathbb Q^n}-\partial_xv_k^{\mathbb Q}\|_\infty
 =\|\partial_xw_k\|_\infty
 \le C\|w_k\|_\infty^{1/2}\|\partial_{xx}w_k\|_\infty^{1/2}\to0
 \qquad\text{as }n\to\infty,
\]
for each $0\le k<m$, proving \eqref{eq:frozengrad}.
\end{proof}

\subsection{Proof of Lemma~\ref{lem:fixedpoint}}\label{app:fixedpoint}
\begin{proof}
Throughout this proof, $C$ denotes a generic constant independent of $\mathbb Q$ and $n$, which may change from line to line. 

We denote by $X$ the canonical process on $C([0,T];\R^d)$ and by $\mathcal K$ the set of probability measures $\mathbb Q$ on $C([0,T];\R^d)$ such that
\[
 \mathbb Q_0=\Law(\xi),\qquad
 \E_{\mathbb Q}\sup_{t\le T}|X_t-X_0|^2\le C_1,\qquad
 \E_{\mathbb Q}\frac{|X_t-X_s|^4}{(1+|X_0|)^4}\le C_2|t-s|^2\quad(s,t\in[0,T]),
\]
where $C_1,C_2>0$ are chosen below. Since $\mathbb Q_0=\Law(\xi)$ and $\xi\in L^2$, we have $\mathcal K\subset\cP_2(C([0,T];\R^d))$. The conditions defining $\mathcal K$ are convex in $\mathbb Q$, so $\mathcal K$ is convex. To see that $\mathcal K$ is weakly closed, let $\mathbb Q^n\in\mathcal K$ converge weakly to $\mathbb Q$. Since $\omega\mapsto\omega_0$ is continuous, $\mathbb Q_0^n\to\mathbb Q_0$ weakly, so $\mathbb Q_0=\Law(\xi)$. For $N>0$, the function $\omega\mapsto\min\{\sup_{t\le T}|\omega_t-\omega_0|^2,N\}$ is bounded and continuous on $C([0,T];\R^d)$, so
\[
 \E_{\mathbb Q}\min\Bigl\{\sup_{t\le T}|X_t-X_0|^2,N\Bigr\}
 =\lim_{n\to\infty}\E_{\mathbb Q^n}\min\Bigl\{\sup_{t\le T}|X_t-X_0|^2,N\Bigr\}
 \le C_1.
\]
Letting $N\to\infty$, monotone convergence gives $\E_{\mathbb Q}\sup_{t\le T}|X_t-X_0|^2\le C_1$. The same argument gives the increment bound, so $\mathbb Q\in\mathcal K$. For tightness, we define the continuous map $\Upsilon(\omega):=(\omega_0,\check\omega)$, with $\check\omega_t:=(\omega_t-\omega_0)/(1+|\omega_0|)$, from $C([0,T];\R^d)$ to $\R^d\times C([0,T];\R^d)$. Its inverse $(x,\check\omega)\mapsto x+(1+|x|)\check\omega$ is continuous on the image of $\Upsilon$. Under $\mathbb Q\in\mathcal K$, $\omega_0$ has the law of $\xi$, $\check\omega_0=0$ and $\E_{\mathbb Q}|\check\omega_t-\check\omega_s|^4\le C_2|t-s|^2$, so the Kolmogorov criterion gives the tightness of the laws of $\Upsilon$ under $\mathbb Q\in\mathcal K$, hence of $\mathcal K$. Thus $\mathcal K$ is weakly compact.

We next show that $\mathcal K$ is nonempty and $\mathcal T\mathcal K\subset\mathcal K$ for $C_1,C_2$ large enough. Let $\mathbb Q\in\mathcal K$ and define $M_t:=\int_0^t\sigma(r,\boldsymbol X_r^{\mathbb Q})\dW_r$. By \eqref{eq:frozenstate}, $X_t^{\mathbb Q}=\xi+\int_0^tb(r,\boldsymbol X_r^{\mathbb Q},u_r^{\mathbb Q})\dd r+M_t$.  By \eqref{eq:uhatbounds}, \eqref{eq:feedbackbounds} and Assumption~\ref{ass:coefficients}\ref{ass:b0}--\ref{ass:drift}, $|b(r,\boldsymbol X_r^{\mathbb Q},u_r^{\mathbb Q})|\le C(1+|\boldsymbol X_r^{\mathbb Q}|)\le C(1+|\xi|)+C\sup_{s\le r}|X_s^{\mathbb Q}-\xi|$. Applying Gronwall's inequality, we obtain
\[
 \sup_{t\le T}|X_t^{\mathbb Q}-\xi|\le C\Bigl(1+|\xi|+\sup_{t\le T}|M_t|\Bigr),\qquad
 \frac{|X_t^{\mathbb Q}-X_s^{\mathbb Q}|}{1+|\xi|}\le C|t-s|\Bigl(1+\sup_{r\le T}|M_r|\Bigr)+|M_t-M_s|.
\]
Since $\sigma$ is bounded, the Burkholder--Davis--Gundy inequality gives
\[
 \E\sup_{t\le T}|X_t^{\mathbb Q}-\xi|^2\le C(1+\E|\xi|^2),\qquad
 \E\frac{|X_t^{\mathbb Q}-X_s^{\mathbb Q}|^4}{(1+|\xi|)^4}\le C|t-s|^2.
\]
The same bounds hold for the solution of $\dd X_t=b(t,\boldsymbol X_t,0)\dt+\sigma(t,\boldsymbol X_t)\dW_t$, $X_0=\xi$. Choosing $C_1,C_2$ large enough, we obtain that $\mathcal K$ is nonempty and $\mathcal T\mathcal K\subset\mathcal K$.

Let $(\mathbb Q^n)_{n\in\mathbb N}\subset\mathcal K$ with $\mathbb Q^n\to\mathbb Q$ weakly. The maps $\omega\mapsto\boldsymbol\omega_t$ and $\omega\mapsto(\omega_{s_1},\ldots,\omega_{s_m})$ are $\sqrt m$-Lipschitz on $C([0,T];\R^d)$. Hence
\[
 \sup_{t\le T}d_{\BL}(\mathbb Q_{\Pi(t)}^n,\mathbb Q_{\Pi(t)})\le \sqrt m\,d_{\BL}(\mathbb Q^n,\mathbb Q),
 \qquad
 d_{\BL}(\mathbb Q_\Pi^n,\mathbb Q_\Pi)\le \sqrt m\,d_{\BL}(\mathbb Q^n,\mathbb Q).
\]
By \eqref{eq:flatcont}, \eqref{eq:frozenhyp} holds, so~\eqref{eq:frozengrad} of Lemma~\ref{lem:frozenproblem} gives $\lim_{n\to\infty}\|\partial_xv_k^{\mathbb Q^n}-\partial_xv_k^{\mathbb Q}\|_\infty=0$ for each $k$. The states $X^{\mathbb Q^n}$ and $X^{\mathbb Q}$ are driven by the same $\xi$ and $W$, and by \eqref{eq:envelope}, the drift of $X^{\mathbb Q}$ is $\partial_yH(t,\boldsymbol X_t^{\mathbb Q},\partial_{x_{k+1}}v_k^{\mathbb Q}(t,\boldsymbol X_t^{\mathbb Q}))$ for $t\in[s_k,s_{k+1})$, and similarly for $X^{\mathbb Q^n}$. By Assumption~\ref{ass:coefficients}\ref{ass:sigma}, $\|\sigma(t,\boldsymbol X_t^{\mathbb Q^n})-\sigma(t,\boldsymbol X_t^{\mathbb Q})\|\le L|\boldsymbol X_t^{\mathbb Q^n}-\boldsymbol X_t^{\mathbb Q}|$. Since $\partial_yH$ is Lipschitz continuous in $(\boldsymbol x,y)$ on $\{|y|\le\bar C\}$ by Assumption~\ref{ass:coefficients}\ref{ass:Hreg}, $|\boldsymbol X_t^{\mathbb Q^n}-\boldsymbol X_t^{\mathbb Q}|^2\le(m+1)\sup_{s\le t}|X_s^{\mathbb Q^n}-X_s^{\mathbb Q}|^2$ and \eqref{eq:feedbackbounds} holds, the Burkholder--Davis--Gundy and Gronwall inequalities give $\lim_{n\to\infty}\E\sup_{t\le T}|X_t^{\mathbb Q^n}-X_t^{\mathbb Q}|^2=0$. Thus $\mathcal T$ is sequentially continuous on $\mathcal K$, hence continuous, since the weak topology on $\mathcal K$ is metrizable by $d_{\BL}$.

We regard probability measures on $C([0,T];\R^d)$ as a subset of the finite signed measures on $C([0,T];\R^d)$, endowed with the locally convex topology induced by integration against bounded continuous functions. By the Schauder--Tychonoff theorem, there exists $\mathbb Q^*\in\mathcal K$ such that
\[
 \mathcal T\mathbb Q^*=\mathbb Q^*.
\]
This proves the desired result.
\end{proof}

\section{Further convex law penalties}\label{fv:sec:examples}

The squared MMD of Section~\ref{sec:softpath} is not the only penalty that satisfies Assumption~\ref{ass:coefficients}\ref{ass:penalties}. We first give a sufficient condition, and then two further examples.

Let $n\ge1$, and let $\mathcal D(\mu,\rho)$ be a discrepancy between a law $\mu\in\cP_2(\R^n)$ and a target law $\rho\in\cP_2(\R^n)$. We denote by $\frac{\delta\mathcal D}{\delta\mu}(\mu,\rho)$ the flat derivative of $\mu\mapsto\mathcal D(\mu,\rho)$, and we assume the following, with $\alpha$ as in Assumption~\ref{ass:coefficients}.
\begin{enumerate}[label=(D\arabic*)]
\item\label{fv:D1} For every $\rho$, the map $\mu\mapsto\mathcal D(\mu,\rho)$ is convex along mixtures.
\item\label{fv:D2} The flat derivative $\frac{\delta\mathcal D}{\delta\mu}(\mu,\rho)$ is bounded in $C_b^{4+\alpha}(\R^n)$, uniformly in $(\mu,\rho)$, and there exists $C\ge0$ such that, for all $\mu,\mu',\rho,\rho'\in\cP_2(\R^n)$,
\begin{equation}\label{eq:examplelawbound}
 \left\|\frac{\delta\mathcal D}{\delta\mu}(\mu,\rho)
       -\frac{\delta\mathcal D}{\delta\mu}(\mu',\rho')\right\|_{C_b^1}
 \le C\bigl(d_{\BL}(\mu,\mu')+d_{\BL}(\rho,\rho')\bigr).
\end{equation}
\item\label{fv:D3} The map $\rho\mapsto\mathcal D(\delta_0,\rho)$ is bounded and continuous with respect to $d_{\BL}$.
\end{enumerate}
Let $\mathcal D$ on $\cP_2(\R^d)$ and $\tilde{\mathcal D}$ on $\cP_2(\R^{md})$ satisfy \ref{fv:D1}--\ref{fv:D3}. We take target laws $\rho_t\in\cP_2(\R^d)$ satisfying~\eqref{eq:targetregularity}, a target law $\rho_\Pi\in\cP_2(\R^{md})$, weights $\lambda_f,\lambda_g\ge0$, and the penalties
\[
 F_t(\mu):=\lambda_f\,\mathcal D(\mu,\rho_t),\qquad
 G(\eta):=\lambda_g\,\tilde{\mathcal D}(\eta,\rho_\Pi).
\]
Then, $F_t$ and $G$ satisfy Assumption~\ref{ass:coefficients}\ref{ass:penalties}. Indeed, they are convex along mixtures by \ref{fv:D1}. Their flat derivatives $f_t(\mu,\cdot)=\lambda_f\frac{\delta\mathcal D}{\delta\mu}(\mu,\rho_t)$ and $g(\eta,\cdot)=\lambda_g\frac{\delta\tilde{\mathcal D}}{\delta\mu}(\eta,\rho_\Pi)$ are bounded in $C_b^{4+\alpha}$ by \ref{fv:D2}, and \ref{fv:D2} together with~\eqref{eq:targetregularity} gives~\eqref{eq:flatcont}. Finally, $F_t(\delta_0)=\lambda_f\mathcal D(\delta_0,\rho_t)$ is bounded and continuous in $t$ by \ref{fv:D3} and~\eqref{eq:targetregularity}.

In each example below, we define $\mathcal D$ and check \ref{fv:D1}--\ref{fv:D3}.
\begin{example}\label{ex:features}
Let $\phi:\R^n\to\R^r$ be in $C_b^{4+\alpha}$, and let $\Psi:\R^r\to\R$ be $ C^2$ and convex, with minimum $\Psi(0)=0$. We define
\[
 \bar\phi_{\mu,\rho}:=\int\phi\,\dd(\mu-\rho),\qquad
 \mathcal D(\mu,\rho):=\Psi(\bar\phi_{\mu,\rho}).
\]
The vector $\bar\phi_{\mu,\rho}$ describes the difference between $\mu$ and $\rho$, and $\Psi$ measures its size. Concrete choices are $\Psi(z)=\sum_{i=1}^r\psi(z_{i})$ with
\[
 \psi(s)=\tfrac12s^2,\qquad
 \psi(s)=\log\cosh s,\qquad
 \psi(s)=\sqrt{\delta^2+s^2}-\delta,\quad \delta>0.
\]
For the first choice $\psi(s)=\frac12s^2$, the discrepancy $\mathcal D(\mu,\rho)=\frac12|\bar\phi_{\mu,\rho}|^2$ is half the squared MMD for the finite-rank kernel $k(x,y)=\phi(x)\cdot\phi(y)$. The other two choices are smooth versions of $|s|$. The coordinates of $\phi$ can be, for example, smooth compactly supported functions.

We now verify \ref{fv:D1}--\ref{fv:D3}. The map $\mu\mapsto \bar\phi_{\mu,\rho}$ is affine and $\Psi$ is convex, so $\mathcal D$ is convex along mixtures \hypersetup{linkcolor=red}, which gives \ref{fv:D1}. Its flat derivative is
\[
 \frac{\delta\mathcal D}{\delta\mu}(\mu,\rho)(x)
 =\nabla\Psi(\bar\phi_{\mu,\rho})\cdot\phi(x).
\]
\hypersetup{linkcolor=red}Since $\phi$ is bounded and Lipschitz, the vectors $\bar\phi_{\mu,\rho}$ lie in a fixed compact set $K\subset\R^r$, and
\[
 |\bar\phi_{\mu,\rho}-\bar\phi_{\mu',\rho'}|
 \le C\bigl(d_{\BL}(\mu,\mu')+d_{\BL}(\rho,\rho')\bigr).
\]
Since $\Psi$ is $C^2$, the functions $\Psi$ and $\nabla\Psi$ are bounded and Lipschitz on $K$. Together with $\phi\in C_b^{4+\alpha}$ and the formula for the flat derivative, this gives \ref{fv:D2}. Since $\mathcal D(\delta_0,\rho)=\Psi(\bar\phi_{\delta_0,\rho})$, it also gives \ref{fv:D3}.
\end{example}

\begin{example}\label{ex:bins}
Let $w_1,\ldots,w_r\in C_b^{4+\alpha}(\R^n)$ be nonnegative with $\sum_{j=1}^r w_j=1$; for instance, $w_1(x)=(1+\tanh(e\cdot x))/2$ and $w_2=1-w_1$ for a fixed $e\in\R^n$. For $0<\eps<1$, we define
\[
 p_j(\mu):=(1-\eps)\int w_j\,\dd\mu+\frac{\eps}{r},\qquad
 \mathcal D_\varphi(\mu,\rho)
 :=\sum_{j=1}^r p_j(\rho)\varphi\left(\frac{p_j(\mu)}{p_j(\rho)}\right).
\]
Thus $\mathcal D_\varphi(\mu,\rho)$ is the $ \varphi$-divergence between the probability vectors $(p_j(\mu))_{j\le r}$ and $(p_j(\rho))_{j\le r}$, obtained from $\mu$ and $\rho$ by smooth binning. The following choices give the KL, Pearson $\chi^2$ and squared Hellinger penalties, respectively~\cite{sreekumar2023}:
\[
\begin{array}{c|c|c}
 &\varphi(s)&\mathcal D_\varphi(\mu,\rho)\\[2pt]\hline
 \mathrm{KL}&s\log s-s+1&\displaystyle\sum_j p_j(\mu)\log\frac{p_j(\mu)}{p_j(\rho)}\\[7pt]
 \chi^2&(s-1)^2&\displaystyle\sum_j\frac{(p_j(\mu)-p_j(\rho))^2}{p_j(\rho)}\\[7pt]
 \mathrm{Hellinger}^2&(\sqrt s-1)^2&\displaystyle\sum_j(\sqrt{p_j(\mu)}-\sqrt{p_j(\rho)})^2
\end{array}
\]
We now verify \ref{fv:D1}--\ref{fv:D3}. Each $p_j$ is affine in $\mu$ and each $\varphi$ above is convex, so these penalties are convex along mixtures, which gives \ref{fv:D1}. Their flat derivatives are
\begin{equation}\label{eq:app flat deri}
 \frac{\delta\mathcal D_\varphi}{\delta\mu}(\mu,\rho)(x)
 =(1-\eps)\sum_{j=1}^r
 \varphi'\left(\frac{p_j(\mu)}{p_j(\rho)}\right)w_j(x).
\end{equation}
Because of the term $\eps/r$, every $p_j(\mu)$ and $p_j(\rho)$ lies in $[\eps/r,1]$. Since $w_j$ is bounded and Lipschitz, $p_j(\mu)$ is Lipschitz in $\mu$ with respect to $d_{\BL}$. Hence the ratios $p_j(\mu)/p_j(\rho)$ lie in the interval $[\eps/r,r/\eps]$, on which $\varphi$, $\varphi'$ and $\varphi''$ are bounded, and they are Lipschitz in $(\mu,\rho)$. Together with $w_j\in C_b^{4+\alpha}$ and \eqref{eq:app flat deri}, we have \ref{fv:D2} and \ref{fv:D3}.
\end{example}

\bibliographystyle{abbrvnat}
\bibliography{ref}
\end{document}